\documentclass[11pt]{article}

\usepackage{amsmath}
\usepackage{amssymb}
\usepackage{latexsym}
\usepackage{amsmath, amsfonts,amssymb, amsthm, euscript,makeidx,color,mathrsfs}
\usepackage{hyperref}
\usepackage{natbib}

\newtheorem{assumption}{Assumption}
\def\qed{ \ \vrule width.2cm height.2cm depth0cm\smallskip}

\newcommand{\la}{\langle}
\newcommand{\ra}{\rangle}
\newcommand{\hP}{\hat\dbP}

\newcommand{\ba}{\begin{array}}
\newcommand{\ea}{\end{array}}
\newcommand{\be}{\begin{equation}}
\newcommand{\ee}{\end{equation}}
\newcommand{\bea}{\begin{eqnarray}}
\newcommand{\eea}{\end{eqnarray}}
\newcommand{\beaa}{\begin{eqnarray*}}
\newcommand{\eeaa}{\end{eqnarray*}}

\def\neg{\negthinspace}

\def\a{\alpha}

\def\g{\gamma}
\def\d{\delta}

\def\m{\mu}
\def\n{\nu}
\def\si{\sigma}

\def\f{\varphi}

\def\o{\omega}

\def\Th{\Theta}

\def\O{\Omega}
\def\cA{{\cal A}}

\def\cD{{\cal D}}

\def\cF{{\cal F}}
\def\cG{{\cal G}}
\def\cH{{\cal H}}

\def\cL{{\cal L}}

\def\sP_1{{\cal P}}

\def\hC{\mathbb{C}}

\def\hE{\mathbb{E}}
\def\hF{\mathbb{F}}
\def\hG{\mathbb{G}}

\def\hL{\mathbb{L}}

\def\hN{\mathbb{N}}

\def\hP{\mathbb{P}}
\def\hQ{\mathbb{Q}}
\def\hR{\mathbb{R}}
\def\hS{\mathbb{S}}

\def\sB{\mathscr{B}}
\def\sC{\mathscr{C}}
\def\sD{\mathscr{D}}

\def\sH{\mathscr{H}}

\def\scL{\mathscr{L}}
\def\scC{\mathscr{C} }
\def\sM{\mathscr{M}}

\def\sP{\mathscr{P}}

\def\sR{\mathscr{R}}

\def\sT{\mathscr{T}}

\def\sX{\mathscr{X}}

\def\no{\noindent}

\def\ss{\smallskip}
\def\ms{\medskip}
\def\bs{\bigskip}
\def\q{\quad}
\def\qq{\qquad}

\def\pa{\partial}
\def\cd{\cdot}
\def\cds{\cdots}
\def\lan{{\langle}}
\def\ran{{\rangle}}

\def\bd{{\bf d}}
\def\bm{{\bf m}}
\def\bn{{\bf n}}
\def\bx{{\bf x}}
\def\by{{\bf y}}
\def\bz{{\bf z}}
\def\bD{{\bf D}}

\def\bp{{\bf P}}
\def\bQ{{\bf Q}}

\def\bZ{{\bf Z}}
\def\vf{{\varphi}}

\def\tr{\hbox{\rm tr}}

\def\qed{ \hfill \vrule width.25cm height.25cm depth0cm\smallskip}
\newcommand{\dfnn}{\stackrel{\triangle}{=}}
\newcommand{\basa}{\begin{assumption}}
\newcommand{\easa}{\end{assumption}}

\newcommand{\bas}{\begin{assum}}
\newcommand{\eas}{\end{assum}}

\def\lan{\mathop{\langle}}
\def\ran{\mathop{\rangle}}

\def\limP2{\,\mathop{\buildrel \Pi_2\over\longrightarrow\,}}

\def\pa{\partial}

 \def\cd{\cdot}
\def\cds{\cdots}

\def\as{\hbox{\rm -a.s.{ }}}

\def\tr{\hbox{\rm tr$\,$}}

\def\ind{{\perp\neg\neg\neg\perp}}

\def\dis{\displaystyle}

\def\bx{{\bf x}}

\def\bP{{\bf P}}
\def\1{{\bf 1}}
\def\by{{\bf y}}

\def\:{\!:\!}
\def\reff#1{{\rm(\ref{#1})}}

at 9pt

\begin{document}

\newtheorem{thm}{Theorem}[section]
\newtheorem{lem}[thm]{Lemma}
\newtheorem{cor}[thm]{Corollary}
\newtheorem{prop}[thm]{Proposition}
\newtheorem{rem}[thm]{Remark}
\newtheorem{eg}[thm]{Example}
\newtheorem{defn}[thm]{Definition}
\newtheorem{assum}[thm]{Assumption}
\newcommand{\ts}{\mathsf{T}}
\renewcommand {\theequation}{\arabic{section}.\arabic{equation}}
\def\thesection{\arabic{section}}

\title{\bf Superpositions for General Conditional  Mckean-Vlasov Stochastic Differential Equations}

\author{
Qi Feng\thanks{\noindent Department of
Mathematics, Florida State University, Tallahassee, 32306; email: qfeng2@fsu.edu. This author is partially supported by the National Science Foundation under grant \#DMS-2420029.}
~ and ~ Jin Ma\thanks{ \noindent Department of
Mathematics, University of Southern California, Los Angeles, 90089;
email: jinma@usc.edu. This author is supported in part
by US NSF grants \#DMS-1908665.} 
}
\date
\maketitle
\date{\today}
\maketitle

\begin{abstract}
In this paper, we study the connection between a general class of {\it Conditional Mckean-Vlasov Stochastic Differential Equations} (CMVSDEs) and its corresponding (infinite dimensional) {\it Conditional Fokker-Planck Equation}. The CMVSDE under consideration is similar to the one studied in \cite{BLM}, which is a non-trivial generalization of the McKean-Vlasov SDE with common noise and is closely related to a new type of non-linear Zakai equation that has not been studied in the literature. The main purpose of this paper is to establish the superposition principles among the three subjects so that their well-posedness can imply each other. More precisely, we shall first prove the superposition principle between the  non-linear Zakai equation,  a non-linear measure-valued stochastic PDE, and a CMVSDE; and then prove the superposition principle between an infinite dimensional  conditional Fokker-Planck equation and the nonlinear Zakai equations. It is worth noting that none of the (weak) well-posedness of these SDEs/SPDEs in such generality have been investigated in the literature. 
\end{abstract}


\vfill \bs

\no


\bs

\no{\it 2000 AMS Mathematics subject classification:} 60H10,15,30;
35R60, 34F05.

\eject

\section{Introduction}
\label{sect-Introduction}
\setcounter{equation}{0}

In this paper we are interested in the so-called {\it superposition principle} among three types of stochastic differential or partial differential equations (SDEs/SPDEs). These SDEs/SPDEs are of the generalized forms of some well-known equations in the literature, but none of them have been systematically studied due to their generalities. In particular, we shall consider a class of 
{\it conditional McKean-Vlasov SDEs} (CMVSDEs), which is closely related to some newly developed stochastic mean-field game/control problems with partial information, the associated measure-valued (nonlinear)  Zakai-type Stochastic PDEs (SPDEs),  and a class of (infinite dimensional) conditional Fokker-Planck equations concerning the law of the aforementioned measure-valued nonlinear Zakai equations and the underlying CMVSDE.  Our main purpose is to establish the  connection among the (weak) well-posedness of these SDEs/SPDEs from the perspective of the  superposition principle.

Let us begin by a brief description of the three types of stochastic differential/partial differential equations that we shall focus on in this paper.

{\bf The Conditional McKean-Vlasov SDE (CMVSDE).} One of the 
main motivations of this paper is the following SDE that can be considered as a generalized form of the (conditional) McKean-Vlasov SDE that often appears in stochastic mean-field control problems, especially under partial information:
\bea
\label{SDE0}
\left\{\ba{lll}
dX_t=b(t, X_t, Y_{\cd\wedge t}, \m^{X|Y}_{\cd\wedge t})dt +\si(t, X_t, Y_{\cd\wedge t}, \m^{X|Y}_{\cd\wedge t})dB^1_t+\rho(t, X_t, Y_{\cd\wedge t}, \m^{X|Y}_{\cd\wedge t})dB^2_t, \ms\\
dY_t= h(t, X_t, Y_{\cd\wedge t}, \m^{X|Y}_{\cd\wedge t})dt+ dB^2_t, \\
X_0=x, Y_0=0.
\ea\right. 
\eea
where $\m^{X|Y}_{t}$, $t\ge 0$, is the conditional law of $X_t$, given $\cF^Y_t\dfnn \si\{Y_s: 0\le s\le t\}$, 
and $(B^1, B^2)$ is a, say,  2-dimensional Brownian motion defined on some probability space $(\O, \cF, \hP; \hF)$. A simpler version of SDE (\ref{SDE0}) has been studied recently in \cite{BLM} (with $\rho\equiv0$), and was named as {\it Conditional McKean-Vlasov SDE} (CMVSDE), which first appeared in the context of mean-field game/control problems with 
common noise when $\rho\neq 0$ but $h\equiv0$ (cf. \cite{carmona2018probabilistic}). 

We should note that in general the well-posedness of the  CMVSDE 
(\ref{SDE0}), even in the weak sense, is a quite subtle issue. In the case of ``common noise" (i.e., $h\equiv 0$) the existence of the solution to CMVSDE is proved in \cite{lacker2020}  via an approximation and limiting argument. In particular, in this case the law 
$\m^{X|Y}$ is conditioned on a fixed Brownian motion, which turns out to be much more technically benign compared to the case $h\neq 0$ as we shall see below. In fact, the 
CMVSDE of the form (\ref{SDE0}) was first proposed in \cite{BLM}, and the weak well-posedness of a slightly simpler form was recently studied in \cite{buckdahn2021general} with additional assumption on function $h$, i.e. $h(t,x,y_{\cdot\wedge t})=\sum_{i=1}^Nf_i(t,x)g_i(,y_{\cdot\wedge t})$, for some bounded measurable functions $f_i\in\hC^{1,2}([0,T],\hR)$ and $g_i$.  
But all these results are a far cry from a complete resolution of the problem. In fact,  it is actually an intriguing problem to find a way to obtain the well-posedness of the CMVSDE (\ref{SDE0}),  and 
this is one of the main purposes  that we pursue the {\it superposition principle} in this paper.

\ss
{\bf The Non-linear Zakai Equation.} If we put CMVSDE (\ref{SDE0}) into a ``nonlinear filtering" framework, then it is conceivable that the conditional law $\m^{X|Y}$ can be determined by a measure-valued SPDE. For the sake of argument let us assume that SDE (\ref{SDE0}) has a (weak) solution $(X, Y)$, defined on some probability space $(\O, \cF, \hP)$, so that the measure-valued process $\m=\m^{X|Y}$ is well-defined on this space. Then, following the standard nonlinear filtering arguments, under appropriate technical conditions, we can define a new probability measure $ \hQ^0\sim \hP$, such that
$(B^1, Y)$ is a $\hQ^0$-Brownian motion, thanks to the Girsanov theorem. In particular, let 
$L_t:=\frac{d\hQ^0}{d\hP}\big|_{\cF_t}$, $t\in[0,T]$, then 
by the well-known {\it Kallianpur-Strieble} formula, it holds that
\bea 
\label{dual}
\lan\m_t, \f\ran =\hE^{\hP}[\f(X_t)|\cF^Y_t]=\frac{\hE^{\hQ^0}[L_t\f(X_t)|\cF^Y_t]}{\hE^{\hQ^0}[L_t|\cF^Y_t]}=\frac{\la \n_t, \f\ra}{\la 
\n_t, 1\ra}, \qq t\in[0,T], \q\f\in \hC_b(\hR),
\eea
where $\n_t(A)=\hE^{\hQ^0}[L_t\1_{\{X_t\in A\}}|\cF^Y_t]$, $A\in \sB(\hR)$, is the so-called {\it unnormalized conditional distribution} of $X_t$ given $\cF^Y_t$. Then, as we shall see in \S3, one can argue that
the measure-valued process $\n$  satisfies the following SPDE under $\hQ^0$:
for $t\in [0,T]$ and  $\f\in \hC^2_0(\hR^d)$,
\bea
\label{Zakai0}
\la \n_t,\f\ra =\la \n_0,\f\ra+\int_0^t\la \n_s,\scL[\f](s, \cd, \cd,\mu) \ra ds+\int_0^t \la \n_s,\sH[\f](s,  \cd, \cd, \mu)\ra dY_s, 
\eea
where, $\scL$ and $\sH$ are differential operators defined respectively by 
\bea
\label{scL}
\left\{\ba{lll}
\dis \scL[\ \!\cd\ \!](t, x, \o, \m):=\frac{1}{2}\big[\tr [(\si\si^\top+\rho\rho^\top)(t, x, \o, \m_{\cd\wedge t})\pa^2_{xx} \big]+( b(t, x, \o, \m_{\cd\wedge t}),\nabla_{x} ),\\
\sH[\ \!\cd\ \!](t, x, \o, \m):=(\rho(t,  x, \o, \m_{\cd\wedge t}),\nabla_x)- h(t, x, \o, \m_{\cd\wedge t}).
 \ea\right.
\eea
In the above 
$(t, x, \o, \m)\in [0,T]\times\hR^d\times \O^0\times \hC([0,T];\sP_1)$, $\sP_1=\sP_1(\hR^d)$ is the space of all probability measures on $\hR^d$, endowed with 1-Wasserstein metric (see \S2 for the detailed definition); $\m_{\cd\wedge t}$ means the paths of $\m$ up to time $t$; and 
$(\cd, \cd)$ denotes the Euclidean inner product. 
We should  note that while the process $\m=\{\m_t\}$ takes values in the space $\sP_1$, the unnormalized conditional distribution process $\n$ is not probability-measure-valued in general. We therefore consider   $\sM_1(\hR^d)$, the space of all finite Borel measures $\n$ on $\hR^d$
with $\int_{\hR^d} |x|\n(dx)<\infty$. Then $\sP_1(\hR^d)\subset \sM_1(\hR^d)$, and in (\ref{Zakai0}) $\la \n,\f\ra:= \int_{ \hR^d }\f(x)\n(dx)$ denotes the {\it dual product} for $\n\in\sM_1(\hR^d)$ and $\f\in\hC_0(\hR^d)$. The solution to SPDE (\ref{Zakai0}) is a  measure-valued processes $\n \in\sM_1(\hR^d)$ such that $\n$ and $\m \in\sP_1(\hR^d)$ appearing on the right hand side of (\ref{Zakai0}) are related by the  Bayes-rule (\ref{dual}). 

It is easy to see that the SPDE (\ref{Zakai0}) is highly  {\it nonlinear} since the 
coefficients contain the measure $\m=\m(\n)$  (by a slight abuse of notation, given the relationship (\ref{dual})). 
In the special case when the coefficients  $b$, $\si$,   $h$ and $\rho$ depend only on $(t,x)\in[0,T]\times\hR^d$,  then the SPDE (\ref{Zakai0}) becomes {\it linear},  namely the usual {\it Zakai equation} that is often seen in the study of nonlinear filtering problems with correlated noises (see, e.g., \cite{bensoussan}). We shall therefore refer to SPDE (\ref{Zakai0}) as the {\it nonlinear Zakai equation} in the rest of the paper, for lack of a better name.

Let us now take a closer look at a popular, but simplified case. That  is, when $h\equiv 0$. In that case  
$Y\equiv B^2$ and $L \equiv 1$,  thus $\n=\m$, and the SPDE (\ref{Zakai0}) is often considered as the limiting version of 
a representative particle in a mean-field interacting system with {\it common noise} $B^2$, as the number of particles tends to infinity. We refer to 
 \cite{coghi2019, crisan2014conditional, hammersley2019weak,lacker2018}
  for some background of such SPDEs. However, from the filtering point of view,  the partial observation case ($h\neq 0$) and the complete observation case ($h\equiv0$) are two fundamentally different cases, and the former is much more technically challenging due to the ``non-Lipschitz" relation (\ref{dual}) between the measure-valued processes $\m=\m^{X|Y}$ and $\n$. 
  
Finally, we would like to point out that the coefficients in (\ref{scL}) and those in (\ref{SDE0}) are actually the same under the so-called {\it canonical} set-up. That is, if we take a {\it canonical} reference probability space $(\O, \cF, \hQ^0)$, where 
$\O^0=\hC([0,T], \hR)$, $\cF=\sB(\O)$,  and $\hQ^0$ the Wiener measure. Let $Y$ be the canonical process: $Y_t(\o)=\o(t)$, $t\in[0,T]$, $\o\in\O$. Then,  $Y$ is a Brownian motion under $\hQ^0$, and we can write  any
 $\hF^Y$-progressively measurable random fields $\phi=b, \si, h, \rho$ as
\bea
\label{eqcoeff}
\phi(t, x, Y_{\cd\wedge t}(\o), \m_{\cd\wedge t})=\phi(t,x,\o_{\cd\wedge t}, \m_{\cd\wedge t})=\phi(t,x,\o, \m_{\cd\wedge t}), \q t\in[0,T],
 \eea 
where $\xi_{\cd\wedge t}$ denotes the paths of process $\xi$ up to time $t$. Clearly, in this sense the structures of the coefficients of (\ref{SDE0}) and (\ref{Zakai0}) are identical. We shall keep such an identification in mind, even when  
 $h\neq 0$, as long as $Y$ is a Brownian motion under $\hQ^0$.

\ms
{\bf The Conditional Fokker-Planck Equation.}  Our second superposition principle concerns the 
SPDE (\ref{Zakai0}) and its corresponding Fokker-Planck equation, which is a stochastic PDE defined on the space of probability laws of 
 $ \sM_1(\hR^d) $-valued solution  $\n=\{\n_t\}$ of \eqref{Zakai0}, which will be referred to as the {\it conditional Fokker-Planck equation} in this paper. 
We shall argue that a solution to the conditional Fokker-Planck equation can be ``lifted" to a weak solution to the SPDE (\ref{Zakai0}) by the superposition principle. For technical reasons, for this part of superposition principle we shall restrict ourself to the case where the coefficients of the SPDE (\ref{Zakai0}) be state-dependent on both $Y$ and $\m$, rather than the general path-dependent form in (\ref{eqcoeff}). 

It is worth noting that the superposition principle for (linear) Zakai equation and the corresponding Fokker-Planck equation on $\sP(\sM_1(\hR^d))$   has been studied only recently (see, e.g., \cite{lacker2020, Qiao2022}).  The main idea is to identify the 
the space $\sM_1(\hR^d)$ to the space
$\hR^\infty(=\hR^\hN)$,  and consider the corresponding Fokker-Planck equation (FPE) as an infinite dimensional PDE (see \S2 for detailed description). Roughly speaking, let $\n:=\{\n_t\}_{t\ge0} $, be an $\ \sM_1(\hR^d)$-valued solution for the Zakai equation (\ref{Zakai0}), and denote its finite dimensional projection on the test functions to be:
$$\sM_1(\hR^d)\ni \n_t\mapsto \la \n_t,\psi \ra:=(\la \n_t,\psi_1\ra,\cdots,\la \n_t,\psi_k\ra ) \in \hR^k ,  \q k\in \hN^+,
$$
where  $ \psi=(\psi_1, \cds, \psi_k) \in (C_c^{\infty}( \hR^d))^k$ is any {\it cylindrical test function}. Let us consider the {\it regular conditional probability distribution} (RCPD) of $\n$ on $(\O^0, \cF^0)$, given $\cF^Y_t$, 
defined as follows, 
 \bea\label{rcpd-y}
 \hP^y_t(\cd):=\hE^{\hQ^0}[\, \cd \, |\cF^Y_t](y)\in \sP(\O^0),\q \hQ^0\textit{-a.e.} \q y\in\O^0,
 \eea
where we denote $\hP^Y_t(\cd):=\hE^{\hQ^0}[\, \cd \, |\cF^Y_t]$. 
 Assuming that $\n \in \hL^2_{\hF^Y}([0,T];\sM(\hR^d))$ and for $\hQ^0$ a.e. $y\in\Omega_0$, define $\bP^Y_t(y):=\bP^y_t$, for $t\in [0,T]$,
and denote $\bP^Y_t:=\hP^Y_t\circ (\n_t)^{-1}\in \sP(\sM(\hR^d))$ to be the conditional probability law of $\n_t$ under $\hP^Y_t$, which we shall  refer to as the {\it unnormalized regular conditional probability distribution} (URCPD) of $X_t$. More precisely, for $t\in[0,T]$,  $\bP^Y_t$ is a $\sP(\sM(\hR^d))$-valued random variable, and for each $A\in \sB(\sM(\hR^d))$, the mapping $(t,y)\mapsto \bP^Y_t(y)(A)$ is $\hF^Y$-progressively measurable.     
We shall argue that, $\{\bP^Y_t:t\ge 0\}$ satisfies the following {\it Conditional Fokker-Planck Equation} on  $\sM(\hR^d)$ in the following sense: for any $k\in\hN$, $f\in\hC^\infty(\hR^k)$,  and $\hQ^0$-a.s., it holds that
{\small 
\bea
\label{FK0}
&&\int_{\sM(\hR^d)}f(\la \bn,\psi\ra)(\bP^Y_t-\bP^Y_0)(d\bn)\\
&=&\int_0^t\int_{\sM(\hR^d)}\Big\{\sum_{i=1}^k\pa_if(\la \bn, \psi\ra) \Big[\la \bn ,\frac{1}{2}\si\si^{\ts}(s, \cd , Y_{\cd\wedge s}, \bm ):\nabla^2 \psi_i(\cd)+(\nabla \psi_i(\cd))^{\mathsf T}b(s, \cd , Y_{\cd\wedge s}, \bm ) \ra \Big]\nonumber \\
&&\qq\qq\q+\frac{1}{2}\sum_{i,j=1}^k\pa_{ij}f(\la \bn, \psi\ra)\Big[\la \bn,(\nabla\psi_i)^{\mathsf T}\rho(s, \cd,Y_{\cd\wedge s}, \bm )+h(s, \cd, Y_{\cd\wedge s}, \bm )\psi_i(\cd)\ra\Big]\nonumber \\
&&\qq\qq\qq\qq\qq\qq\times \Big[\la \bn,(\nabla\psi_j)^{\mathsf T}\rho(s, \cd,Y_{\cd\wedge s},  \bm )+h(s, \cd, Y_{\cd\wedge s}, \bm )\psi_j(\cd)\ra\Big]
\Big\}\bP^Y_s(d\bn) ds\nonumber\\
&&+\int_0^t\int_{\sM(\hR^d)} \Big\{\sum_{i=1}^k\pa_i f(\la \bn, \psi\ra)\Big[\la \bn,(\nabla \psi_i(\cd))^{\mathsf T}\rho(s, \cd,Y_{\cd\wedge s}, \bm)+h(s, \cd, Y_{\cd\wedge s}, \bm)\psi_i(\cd)\ra\Big]  \Big\}\bP^Y_s(d\bn)dY_s.\nonumber
\eea}
Here in (\ref{FK0}) $\bm=\bm(\bn)\in\sP(\hR^d)$ is  defined by the Kallianpur-Strieble formula $\la \bm,\f\ra =\frac{\la \bn,\f \ra}{\la \bn,\f\ra}$ for test function $\f\in C_c^{\infty}(\hR^d)$.
To the best of our knowledge, the superposition principle between the (conditional) Fokker-Planck equation (\ref{FK0})  and the 
nonlinear Zakai equation (\ref{Zakai0}) 
in such general form 
is new. In particular, our results will be the non-trivial extensions of the ``common noise" case (i.e., $h\equiv0$) studied in recent  works \cite{barbu2020nonlinear,lacker2020, Qiao2022,  ren2020space}.

The rest of paper is organized as follows. In section 2 we introduce some necessary notions and definitions, along with some basic facts about 
CMVSDE, nonlinear filtering, and $L$-derivatives. In section \ref{CMVSDE to Zakai}, we show the existence of non-linear Zakai equation \eqref{Zakai0} given the solution of the conditional McKean-Vlasov SDE. In section \ref{relation 2}, we show the existence of the weak solution of the conditional McKean-Vlasov SDE given the solution of the non-linear Zakai equation. 
 In section \ref{relation 1}, we derive the conditional Fokker-Planck equation for the non-linear Zakai equation where we reduce to state dependent coefficients for the non-linear Zakai equation. In the end, we prove the superposition principle for the non-linear Zakai equation.

\section{Preliminaries}
\setcounter{equation}{0}
Throughout this paper we denote $\hC_T^d:=\hC([0,T]; \hR^d)$, $d\in\hN$, the $\hR^d$-valued continuous functions on $[0,T]$, equipped with the usual sup-norm. (The case $d=\infty$ will be discussed separately below.)

Consider the canonical probability space $\O^0:=\hC_T^{d+1}$, $\cF^0:=\sB(\hC^{d+1}_T)$,  $\hF^0=\{\cF^0_t\}_{t\in[0,T]}:=
\{\sB_t(\hC_T^{d+1})\}_{t\in[0,T]}$, where $\sB_t(\hC^{d+1}_T):=\si\{\o(\cd\wedge t): \o\in\hC^{d+1}_T\}$. We denote $\hP^0_d$ to be the Wiener measure on $\hC_T^d$, and $\hQ^0=\hP^0_d\otimes\hP^0_1$, the Wiener measure on $\hC^{d+1}_T$. For the sake of our discussion, in what follows we shall denote the {\it canonical process} by $(B^1, Y)$, and denote the generic element of $\O^0$ by 
$\o=(\o^1, y)\in \hC^d_T\times \hC_T=\hC^{d+1}_T$. Then, $(B^1_t, Y_t)(\o)=(\o^1(t), y(t))$, $\o=(\o^1, y)\in \hC^{d+1}_T$, $t\in [0,T]$. Clearly, under $\hQ^0$, $(B^1, Y)$ is a $(d+1)$-dimensional-Brownian motion. We shall also assume that the filtration $\hF^0$ is augmented by all the $\hQ^0$-null sets of $\O^0$ so that it satisfies the {\it usual hypotheses} (cf. e.g., \cite{protter}).  

Now for any sub-$\si$-field $\cG
\subseteq \cF^0$,  and any generic metric space $(\sX, d)$ with $\sB(\sX)$ denoting its topological Borel field, we denote  $\hL^p_{\cG}(\sX)$, $p\ge 1$, to be the space of all $\sX$-valued, 
$\hL^p$-integrable $\cG$-measurable random variables. Similarly, for any sub-filtration $\hG\subseteq\hF$, and $p\ge 1$, we denote
 $\hL^p_\hG([0,T];\sX)$ the space of all $\sX$-valued, $\hG$-adapted processes defined
on $[0,T]$. 
Next, let $\sP(\sX)$ be the space of probability measures on the metric space $\sX$, and let 
 $\sP_1(\sX)=\{\g\in\sP(\sX): \int_{\hR^d}|z|\g(dz)<+\infty\}$,  endowed with the
 {\it 1-Wasserstein metric}:
\bea
\label{Wass1}
W_1(\g_1,\g_2)\neg:=\neg\inf\neg\Big\{\neg\int_{\sX^2}\neg|z_1-z_2|\rho(dz_1 dz_2) :\rho\in\sP(\sX^2),\rho(\cdot\neg\times\neg\sX)=\g_1,\rho(\sX\neg\times\neg\cdot)=\g_2\Big\}.
\eea
It is well-known that $(\sP_1(\sX),W_1)$ is a complete metric space. 
 In particular, if $\sX=\O^0=\hC_T^{d+1}$, we shall simply denote $\sP_1=\sP_1(\O^0)$. Clearly, in this case we have, for $\g_1, \g_2\in\sP_1$, 
\bea
\label{W1}
W_1(\g_1,\g_2)=\inf\{\hE^{\hQ^0}[|\xi_1-\xi_2|]: \xi_i\in \hL^1(\O^0;\hQ^0), \ \hQ^0\circ\xi_i^{-1}=\g_i, ~i=1,2\}. 
\eea

For  fixed $t\in[0,T]$, let us denote the {\it regular conditional probability} $\hQ^{y}_t(A):=\hE^{\hQ^0}\{\1_A|\cF^Y_t\}(y)$, for $A\in \cF^0$, and $\hQ^Y_t(y)=\hQ^y_t$, for $\hQ^0$-a.e. $y\in\Omega_0$. 
Then the mapping $y\mapsto \hQ^{y}_t(A)$ is $\sB_t(\hC_T)/\sB(\hR)$ measurable for any $A\in\cF^0$. Since $\sB_t(\hC_T)$ is generated by the paths $Y_{\cd\wedge t}(y)=y(\cd\wedge t)$, $y\in\hC_T$, we have $\hQ^{y}_t=\hQ^{Y}_t(y)=\hQ^{Y}_t( y_{\cd\wedge t} )$, when there is no confusion. For any   $\xi\in \hL^1_{\cF^0}(\O^0)$ and 
$t\in[0,T]$ we denote:
\bea
\label{PY}
\hP^{Y_{\cd\wedge t}}_\xi (y) (\, \cd \,)=\hP^{y_{\cd\wedge t}}_\xi(\, \cd \,)=\hQ^{y}_t\circ \xi^{-1}(\cd)=\hQ^0[\xi\in \cd\, |\cF^Y_t](y)\in \sP_1,
\eea
which exists $\hQ^0$-a.s.  $t\neg\in\neg[0,T]$. We now consider the mapping $(t, y)\mapsto \hP^{y_{\cd\wedge t}}_\xi
(\cd)\in \sP_1$. First, it was shown in \cite{BLM} that 
the mapping $y\mapsto \hP^{Y_{\cd\wedge t}}_\xi(\cd)(y)=\hP^{y_{\cd\wedge t}}_\xi(\cd)$
is $\sB_t(\hC_T)/\sB(\sP_1)$-measurable, for  $t\in[0,T]$.
On the other hand, for  $y\in \hC_T$, the mapping
$t\mapsto \hP^{y_{\cd\wedge t}}_\xi$ is a $\sP_1$-valued continuous function, which then renders the joint measurability of the mapping $(t, \o^2)\mapsto \hP^{y_{\cd\wedge t}}(\cd)$. 
{In what follows we shall also focus on the space $\sC_T(\sM_1)$,  where $\sM_1=\sM_1(\hR^d)$ is the space of all $\sM_1$-valued continuous functions defined 
on $[0,T]$, and thus $\sP_1\subset \sM_1$. For a given sub-filtration $\hG\subset \hF$, we denote
$\hL^0_\hG(\sC_T(\sM_1))$ to be the space of all $\sM_1$-valued, $\hG$-adapted
continuous processes; and 
\bea
\label{SpG}
\hS_{\hG}^p( \sM_1 ):=\{z\in \hL^0_{\hG}(\sC_T(\sM_1)): 
\hE^{\hQ^0}\Big[\sup\limits_{t\in[0,T]}W_1(z_t,\hQ^0)^p\Big]<+\infty\}.
\eea
The space  $\hS^2_{\hF^Y}( \sM_1)$ (as well as $\hS^2_{\hF^Y}(\sP_1)$) 
will  be particularly useful in our discussion.}

{ \ms
{\bf $L$-derivatives.} We next introduce the notion of the so-called {\it $L$-derivative}, first introduced in \cite{cardaliaguet2010notes} (see \cite{carmona2018probabilistic} for more details), which plays essential role for discussing the space of distributions on the measure-valued processes and allows us to properly define the conditional Fokker-Planck equation \eqref{FK0} on $\sM_1$. We note that in the case of conditional Fokker-Planck equation for the Zakai-type of SPDEs, the underlying dynamics will be the unnormalized conditional law, hence (finite) measure-valued processes. The corresponding $L$-derivatives will be slightly modified from those introduced recently in \cite{lacker2020, Qiao2022,ren2020space}, which mainly deals with the probability measures.  

\begin{defn}
\label{L derivative}
We define $\hC_b^2(\sM_1)$ to be the set of all bounded, continuous functions $F: \sM_1\rightarrow \hR$ such that:

{\rm (i)}  There exists a bounded function $\pa_{\m}F:\sM_1(\hR^d)\times\hR^d\rightarrow \hR$ such that 
\beaa 
\lim_{\varepsilon\rightarrow 0}\frac{F(\m+\varepsilon \xi)-F(\m)}{\varepsilon}=\int_{\hR^d}\pa_{\m}F(\m,x)\xi(dx),  \qq \m,\xi\in \sM_1(\hR^d). \eeaa

\ss
{\rm (ii)} The function $x\mapsto \pa_{\m}F(\m,x)$ is continuously differentiable for $x\in\hR^d$, and the gradient $\nabla_x \pa_{\m}F(\m,x)$ is uniformly bounded in $(\m,x)$, denoted as $\bD_\m F(\m,x)$.

\ss
{\rm (iii)} The second order spacial derivative is defined as $\bD_{\m}^2 F(\m,x,x')$, and  is bounded and continuous for $x,x'\in\hR^d$.

\ss
{\rm {(iv)}} The mixed derivative $\bD_x\bD_{\m}F(\m,x)$ is defined as gradient of $\bD_{\m}F(\m,x)$ and assume that it is continuous and bounded in $(\m,x)$.
\end{defn}

In particular, we shall often consider the $k$-{\it cylindrical} test function $F\in\hC^2_b(\sM_1)$ for some $k\in\hN$, $\f_1, \cds, \f_k\in\hC(\hR^d)$, and $f\in\hC^2(\hR^k)$, such that  we have
\bea 
F(\m)=f(\la \f, \m\ra)=f(\la\m,\f_1\ra, \cdots, \la\m,\f_k\ra ), \qq \m\in\sM_1(\hR^d).
\eea 
By Definition \ref{L derivative}, with a slight  abuse of notation, in this case one can check
\bea
\label{Lcyl}
\bd_{\m}F(\m,x)&=&\sum_{i=1}^k\pa_if(\la\m,\f\ra)D \f_i(x),\quad 
\bd_{x}\bd_{\m}F(\m,x)=\sum_{i=1}^k\pa_if(\la\m,\f\ra)D^2\f_i(x),\nonumber \\
\bd_{\m}^2F(\m,x,x')&=&\sum_{i,j=1}^k\pa_{ij}f(\la\m,\f\ra)D\f_i(x)D\f_j(x')^\top,\quad D\f_i =(
	\nabla\f_i, 	\f_i)^\top.
\eea

Then, the conditional Fokker-Planck equation \eqref{FK0} is equivalent to the following form by using the $L$-derivative and a $k$-cylindrical test function, for $\hQ^0$-a.s.,
\bea\label{FP L}
\int_{\sM(\hR^d)}F (\bn)(\bP^Y_t-\bP^Y_0)(d\bn)
&=&\int_0^t\int_{\sM(\hR^d)} \cL_s^{Y, \bm} F(\bn) \bP^Y_s(d\bn) ds\nonumber \\
&&+\int_0^t\int_{\sM(\hR^d)}\bd_{\bn}F(\bn,x)\cdot \cH(s, x , Y_{\cdot\wedge s},{\bm}) \bP^Y_s(d\bn) dY_s ,
\eea
where the operator $\cL^{Y, \bm}_t$, for $t\in [0,T]$, is defined as below, 
\bea \label{FP operator}
\cL_t^{Y, \bm} F(\bn)&:=&\int_{\hR^d}\Big[\bd_{\bn}F(\bn,x)\cdot\bar b(t, x , Y_{\cdot\wedge t},{\bm})+\frac{1}{2}\bd_x\bd_{\bn}F(\bn,x):\bar\si\bar \si^{\ts}(t,x , Y_{\cdot\wedge t},{\bm})\Big]\bn(dx)\nonumber \\
&&+\frac{1}{2}\int_{\hR^d}\int_{\hR^d}\bd_{\bn}^2F(\bn,x,x'):\cH(t,x, Y_{\cdot\wedge t}, \bm)\cH(t,x',Y_{\cdot\wedge t},\bm)^{\mathsf T}\bn(dx)\bn(dx'),
\eea 
and we denote 
\beaa 
\bar b=\begin{pmatrix}
	b\\
	0
\end{pmatrix},\quad 
\bar\sigma=\begin{pmatrix}
	\sigma\\
	0
\end{pmatrix},\quad 
\cH=\begin{pmatrix}
	\rho\\
	h
\end{pmatrix}.                                                                                                                                                                                                                                                                                                                                                                                                                                                                                                                                                                                                                                                                                                                                                                                                                                                                                                                                                                                                                                                                                                                                                                                                                                                                                                                                                                                                                                                                                                                                                                                                                                                                                                                                                                                                                                                                                              
\eeaa

\ms
{\bf Nonlinear Filtering Revisited.} We now briefly recall some facts of nonlinear filtering, as it  motives the main subjects discussed in this paper. By a {\it signal-observation pair  with correlated noise} we mean the solution $(X, Y)$ to the following SDEs defined on  $(\O^0, \cF^0, \hP)$,  $\hP\in\sP(\O^0)$:
\bea
\label{Y0}
\left\{\ba{lll}
dX_t=b(t, X_t)dt+\si(t, X_t)dB^1_t +\rho(t, X_t)dB^2_t,\\
dY_t= h(t,X_t)dt+dB^2_t, 
\ea\right.
\qq t\in[0, T],
\eea
where $(B^1, B^2)$ is a $(\hP, \hF)$-Brownian motion.
Assume that the function $h$ is, say,  bounded,  and let $\bar{L}$ be the solution to the linear SDE:
\bea
\label{L}
d\bar{L}_t= h(t, X_t)\bar{L}_tdB^2_t, \qq t\in[0,T]; \q \bar{L}_0=1, \q \hP\as\neg\neg,
\eea
then one can define a new measure by $\frac{d\hQ^0}{d\hP}\big|_{\cF_T}=\bar{L}_T$, so that $\hP\sim\hQ^0\in\sP(\O^0) $. Futhermore, if we denote $\n_t$ to be the so-called {\it unnormalized conditional law}  in the sense that $\lan\n_t, \f\ran=\hE^{ \hQ_0}[\bar{L}_t\f(X_{ t})|\cF_t^Y]$, $t\in[0,T]$, for any $\f\in \hC_0(\hR^d)$; and $\m_t=\hP^{X_t|Y}$, the conditional law (under $\hP$) of $X_t$, given 
$\cF^Y_t:=\si\{Y_s:0\le s\le t\}$. Then,  $\n$ satisfies the Zakai equation (\ref{Zakai0}) when $b, \si, \rho, h$ depend only on 
$(t, x)$, and $\n$ and $\m$ are exactly related by (\ref{dual}), the well-known {\it Kallianpur-Strieble formula} (or the Bayes rule, see \cite{bensoussan}).

\ms
{\bf Diffusions on $\hR^\infty$.} To end this section, we introduce the space $\hR^\infty$, and establish a framework on which the superposition principle (or the conditional Fokker-Planck equation (\ref{FK0})) of the measure-valued SPDE can be discussed  as a diffusion on $\hR^\infty$. Most of the definitions and results here can be found in, e.g., \cite{AGSbook, Qiao2022}.

To begin with, let us endow $\hR^\infty=\hR^{\hN}$ with the product topology  and the following metric:
\bea
\label{Rinf}
d_{\infty}(x, y):=\sum_{k=1}^\infty\frac{|x^k-y^k|\wedge 1}{2^k}, \q x=(x^k), y=(y^k)\in\hR^\infty.
\eea
Next, we consider the space $\hC_T^ \infty:=\hC([0,T];\hR^\infty)$, as the extension of $\hC^d_T$, with the metric: 
\bea
\label{Cinf}
d_{C, \infty}(\bx, \by):=\sum_{i=1}^\infty\frac{\sup_{t\in[0,T]}(|\bx^k_t-\by^k_t|\wedge 1)}{2^k}, \q \bx=(\bx^k), \by=(\by^k)\in\hC_T^\infty.
\eea
Then, it is known (cf. \cite{AGSbook}) that both $(\hR^\infty, d_{\infty})$ and $(\hC_{T,\infty}, d_{C, \infty})$ are complete, separable metric spaces. 
Next, for $n\in\hN$, we define the canonical projections from $\hR^\infty$ to $\hR^n$ by $\pi^n(x)=(x^1, \cds, x^n)$, $x\in\hR^\infty$. We say that a function $F: \hR^\infty\mapsto \hR$ is a {\it cylindrical function} if there exists $n\in\hN$, such that
$$ F(x)=F(\pi^n(x))=f( \pi^n(x)), \qq x\in[0,T]\times \hR^\infty, 
$$
for some $f\in\hC^n$, and in this case we say that $F$ is {\it $n$-cylindrical}. We denote $\sC_{ n}$ to be the set of all $n$-cylindrical functions, and $\sC_{ \infty}=\bigcup_{n=1}^\infty\sC_{ n}\subset \hC_{  \infty}$  the set of all cylindrical functions. 
Further, for $k, d\in\hN$, we denote $\hC^{(k)}_d=\hC^{(k)}(\hR^d)$ and $\sC^{(k)}_d\subset \hC^{(k)}(\hR^\infty)$ to be the $k$-th continuously differentiable $n$-cylindrical functions.
Finally, we define $\sC^{(k)}_\infty:=\bigcup_{d=1}^\infty\sC^{(k)}_d$.

Now, let us fix a countable dense set $\{\vf_\ell, \ell\in\hN\}\subset \hC_c^{\infty}(\hR^d)=\{\f \in \cap_{k=1}^{\infty} \hC^{(k)}(\hR^d): \mbox{$\f$  has compact support}\}$, and consider the mapping
\bea
\label{sT}
\sT:\sM(\hR^d)\mapsto 
\hR^\infty, \qq \sT(\m)=(\lan\m, \vf_1\ran, \cds, \lan\m, \vf_\ell\ran, \cds). 
\eea
Then clearly,  $\sT$ is a linear mapping, and is injective. Next, define a metric on $\sM(\hR^d)$ by
\bea
\label{dM}
{\bf d}(\m, \n):=\sum_{k=1}^\infty\frac{|\lan \m, \vf_k\ran-\lan \n, \vf_k\ran|\wedge 1}{2^k}, \qq \m, \n\in\sM(\hR^d), 
\eea
then the metric $\bf d$ induces the weak convergence topology on $\sM(\hR^d)$, and thus $\sT(\sM(\hR^d))$ is a closed subset of $\hR^\infty$. Since $\sT$ is linear, we can identify  $\sT(\sM(\hR^d))$ as a closed subspace  of $\hR^\infty$, denoted by $\hR^\infty_{\sM}$. Then $\sT^{-1}: \hR^\infty_\sM\mapsto \sM(\hR^d)$ is well-defined and continuous. In other words, the mapping $\sT$ is an isomorphism between the spaces $(\sM(\hR^d), {\bf d})$ and $(\hR^\infty_{\sM}, d)$.

Next, given $F\in\sC_{T, n}$, we can write $F(x)=F(\pi^n(x))=f(x)$, for some $f\in\hC(\hR^n) $. Thus, for any $\m\in\sM(\hR^d)$, we have
$$ F\circ \sT(\m):=F(\sT(\m))=F(\lan \m, \vf_1\ran, \cds, \lan\m, \vf_n\ran, \cds)=f(\lan\m, \vf_1\ran, \cds, \lan\m, \vf_n\ran).  
$$
{Consequently, by using the mapping $\sT$ and its inverse $\sT^{-1}$  we are able to transform the conditional Fokker-Planck equation \eqref{FK0} into the following equivalent form on $\hR^{\infty}$. More precisely, for $t\in[0,T]$, note that the solution for \eqref{FK0} $\bP_t^Y$ with  $\bP^Y_t(y)=\bP_t^y=\hP^y_t\circ (\n_t)^{-1}\in \sP(\sM(\hR^d))$, $y\in\O^0$, where $\n_t\in \sM(\hR^d)$, would naturally induce a conditional probability law on $\sP(\hR^\infty_\sM)\subseteq \sP(\hR^\infty)$ via the mapping
$\sT$, denoted by $\mathbf Q_t^Y$. That is, $\bQ_t^Y(y)=\bQ^y_t=\bP^y_t\circ \sT^{-1}=\hP^y_t\circ (\n_t)^{-1}\circ\sT^{-1}\in\sP(\hR^\infty)$. Now for any $\m\in\sM(\hR^d)$ and $\varphi=(\varphi_\ell)_{\ell\in\hN}\in\sC_{T,\infty}$, 
we have $\bz=\{\la \m,\varphi_i\ra \}_{i=1}^\infty\in\hR^\infty$. We can now write down the counterpart of conditional Fokker-Planck equation on $\hR^\infty$: for any $k\in \hN$ and $f\in\hC(\hR^k) $, for $\hQ^0$-a.s., 
\bea
 \label{R infinity FP}
\int_{\hR^{\infty}} f(\bz)(\mathbf Q^Y_t -\mathbf Q^Y_0)(d\bz)&=&\int_0^t\int_{\hR^{\infty }} 
\Big\{\sum_{i=1}^{k}\pa_if(\bz)\beta_i(s,Y_{\cdot\wedge s},\bz) \Big\}\mathbf Q^Y_s(d\bz)ds\\
&&+\int_0^t\int_{\hR^{\infty}}\frac{1}{2}\sum_{i,j=1}^k\pa_{ij}f(\bz)\alpha_{ij}(s,Y_{\cdot\wedge s},\bz)\mathbf Q^Y_s(d\bz) ds \nonumber \\
&&+\int_0^t\int_{\hR^{\infty}} \Big\{\sum_{i=1}^k\pa_i f(\bz)\gamma_i(s,Y_{\cdot\wedge s},\bz) \Big\}\mathbf Q^Y_s(d\bz)dY_s, \q t\in[0,T]\nonumber.
\eea 
In the above, the coefficients $\a_{ij}$'s, $\beta_i$'s and $\gamma_i$'s  are defined by 
{\small
\bea \label{alpha beta}
\left\{\ba{lll}
\alpha_{ij}(s, Y_{\cdot\wedge s} ,\bz)  := \Big[\la \sT^{-1}(\bz),(\nabla\f_i(\cdot))^\top\rho(s, \cd,Y_{\cdot\wedge s}, \bm(\sT^{-1}(\bz)))+h(s, \cd, Y_{\cdot\wedge s},\bm(\sT^{-1}(\bz)))\f_i(\cd)\ra\Big]\\
 \qq\qq\qq \times \Big[\la \sT^{-1}(\bz),(\nabla\f_j(\cdot))^\top\rho(s, \cd,Y_{\cdot\wedge s}, \bm(\sT^{-1}(\bz)))+h(s, \cd, Y_{\cdot\wedge s}, \bm(\sT^{-1}(\bz)))\f_j(\cd)\ra\Big],\\ 
\beta_i(s,Y_{\cdot\wedge s},\bz)  := \la \sT^{-1}(\bz) ,L^Y_{s,\bz}\f_i \ra,\\
\gamma_i(s,Y_{\cdot\wedge s},\bz)  := \la \sT^{-1}(\bz),(\nabla \f_i(\cd))^\top\rho(t, \cd,Y_{\cdot\wedge s}, \bm(\sT^{-1}(\bz)))+h(t, \cd, Y_{\cdot\wedge s}, \bm(\sT^{-1}(\bz)))\f_i(\cd)\ra,
\ea\right.  
\eea 
}
where 
$$L^Y_{s,\bz}:=\frac{1}{2}\si\si^{\top}(s, \cd , Y_{\cdot\wedge s},\bm(\sT^{-1}(\bz) )):\nabla^2 \f_i(\cd)+(\nabla \f_i(\cd))^\top b(s, \cd , Y_{\cdot\wedge s}, \bm(\sT^{-1}(\bz) )),$$
and   the mapping $\bm\circ\sT^{-1}: \hR^{\infty}\rightarrow \sP_1(\hR^d)$  is defined by 
$\la \bm(\sT^{-1}(\bz)),\f\ra := \frac{\la \sT^{-1}(\bz),\f\ra  }{\la \sT^{-1}(\bz),1\ra }$, 
 $\f \in\hC^d_T$.
 }

\section{CMVSDE and its nonlinear Zakai equation}\label{CMVSDE to Zakai}
\setcounter{equation}{0}

In this section we derive a Zakai-type of SPDE associated to the  
{\it conditional McKean-Vlasov SDE} (CMVSDE) (\ref{SDE0}). Such a parity resembles the well-known relationship between the nonlinear filtering problem (\ref{Y0}) and the corresponding Zakai equation (\ref{Zakai0}), and will be interpreted as the first ``superposition principle" that we will pursue in the next section. To be more precise, let us recall the SDE syetem (\ref{SDE0}), and note that it is slightly more general than the one studied in \cite{BLM}:
\bea
\label{SDE1}
\left\{\ba{lll}
dX_t= b(t,  X_t, Y_{\cd\wedge t}, \m_{\cd\wedge t})dt+  \si(t,  X_t, Y_{\cd\wedge t}, \m_{\cd\wedge t})dB^1_t+\rho(t,  X_t, Y_{\cd\wedge t}, \m_{\cd\wedge t})dB^2_t;\ms\\
dY_t= h(t,  X_t, Y_{\cd\wedge t}, \m_{\cd\wedge t})dt + dB^2_t, \q\qq  X_0=x, \q Y_0=0,
\ea\right.
\eea
defined on some filtered probability space $(\O, \cF, \hP;\hF=\{\cF_t\})$, where $\m_t =\hP^{X|Y}_t$ is the conditional law of $X_t$ given $\cF^Y_t:=\si\{Y_s:0\le s\le t\}$ under $\hP$,   and
$(B^1, B^2)$ is a $(\hP, \hF)$-Brownian motion. We shall make use of the following standing assumptions on the coefficients:
\begin{assum}
\label{assump1}
(i) The function $\f=(b,\si, h, \rho):[0,T]\times \hR^d\times \hC^d_T\times \sC_T(\sP_1)\mapsto  \hR^{d\times 1}\times \hR^{d\times d}\times \hR\times \hR^{d\times 1}$ is bounded, measurable, and  for some constant $C>0$, it holds 
that
\bea
\label{Lip}
&&|\f(t, x, \by, \m)-\f(t, x', \by, \m')|x
\le C\Big[ |x -x' | + W_1(\m_t, \m'_t)\Big], \\
&&\qq\qq\qq\qq\qq\qq\qq t\in[0,T], ~x, x'\in \hR, ~\by\in \hC_T, ~\m,\m'\in\sC_T(\sP_1);\nonumber
\eea

(ii) For fixed $(t, \by, \m)\in [0,T]\times\hC_T\times \sC_T(\sP_1)$, $\si(t, \cd, \by, \m)\in  \hC^{(2)}_b(\hR^d)$;

\ms
(iii) For some $\si_0>0$, it holds that 
\bea
\label{nondeg}
\si\si^\top(t,x,\by, \m) - \rho\rho^\top(t,x,\by, \m)\ge \si_0, \q (t,x,\by, \m)\in[0,T]\times\hR^d\times \hC^d_T\times\sC_T(\sP_1).
\eea
\end{assum} 

In light of the result in \cite{BLM}, it is a tall order to find the strong solution to a CMVSDE, we shall thus consider the weak solution of 
(\ref{SDE1}). We first recall the definition (see \cite{BLM}):
\begin{defn}
An eight-tuple $(\O, \cF, \hP,\hF,  B^1, B^2, X, Y)$ is called a weak solution to the CMVSDE (\ref{SDE1}) if 

{\rm (i)}  $(\O,\cF,\hP;\hF)$ is a filtered probability space satisfying the usual hypothese;

\ss
{\rm (ii)} $(B^1,B^2)$ is an $(\hP, \hF)$-Brownian motion;

\ss
{\rm (iii)} SDE (\ref{SDE1}) holds with $\m_t=\hP^{X|Y}_t$, for all  $t\in[0,T]$, $\hP$-a.s.
\qed
\end{defn}

It is often useful to recall the usual scheme in seeking a weak solution. Let us consider
the following  SDEs on the canonical space $(\hQ^0, \cF^0, \hQ^0)$ so that $(B^1, Y)$ is a $\hQ^0$-Brownian motion:
\bea
\label{SDE2}
\left\{\ba{lll}
dX_t=\tilde{b}(t,  X_t, Y_{\cd\wedge t}, \m_{\cd\wedge t})dt+ \si(t,  X_t, Y_{\cd\wedge t}, \m_{\cd\wedge t})dB^1_t+\rho(t, X_t, Y_{\cd\wedge t}, \m_{\cd\wedge t})dY_t;\ms\\
dL_t= h(t,  X_t, Y_{\cd\wedge t}, \m_{\cd\wedge t})L_tdY_t,  \qq \q  X_0=x, \q L_0=1.
\ea\right.
\eea
where $\tilde b=b-\rho h$, and $\m\in \hS^2_{\hF^Y}(\sP_1)$.
Clearly, under Assumption \ref{assump1}, given any $\m\in \hS_{\hF^Y}^2(\sP_1)$, the SDE (\ref{SDE2}) has a unique strong solution, denote it by $(X^\m, L^\m)$, on the probability space $(\O, \cF, \hF^{B^1, Y}, \hQ^0)$.  
We can then define a new probability measure $\bar\hP$ on $(\O^0, \cF^0)$ by: $\frac{d\bar\hP}{d\hQ^0}|_{\cF_T}:=L_T $, and define
a regular conditional probability measure via the Kallianpur-Striebel (Bayes) formula: for $t\in[0,T]$, 
\bea
\label{mubar}
\bar\m_t(A)\dfnn \bar\hP\{X^\m_t\in A|\cF^Y_t\}=\frac{\hE^{\hQ^0}\big[L^{\m}_t{\bf 1}_{\{X^\m_t\in A\}} |\cF^Y_t\big]}
{\hE^{\hQ^0}\big[L^{\m}_t |\cF^Y_t\big]}, \qq A\in \sB(\O).
\eea
Note that the process $(Y, \bar \m)$ are $\hC_T\times \sC_T(\sP_1)$-valued and  is independent of $B^1$, under $ \hQ^0$. 
Furthermore, if it turns out that $\m_t=\bar\m_t$, $t\in[0,T]$, $\hQ^0$-a.s., where $\bar\m$ is given by (\ref{mubar}), then by defining 
$B^2_t=Y_t-\int_0^t h(s, X_s, Y_{\cd\wedge s}, \m_{\cd\wedge s})ds$, we see that $(B^1, B^2)$ is a $\bar\hP$-Brownian motion, 
and $\m_t=\bar\hP^{X|Y}_t$, $t\in[0,T]$, is the conditional law of $X=X^\m$ given $\cF^Y_t$, under $\bar\hP$. That is, $(X^\m, Y, \m)$ is a solution to 
(\ref{SDE0}) under $\bar\hP$. We shall therefore call the solution to (\ref{SDE2}),  $(X^\m, L^\m, \m)\in \hL^2_{\hF}(\O^0; \hC_T\times\sC_T(\sP_1))$, a ``canonical solution" to a CMVSDE (\ref{SDE0}), which will be useful in some of our discussions. 

We should note that   the (weak) well-posedness of the CMVSDE (\ref{SDE1}) (that is, finding the fixed point $\m=\bar\m$ in the above heuristic argument) is by no means trivial, we refer to \cite{BLM}  for detailed discussion. But as we shall see in the next section,  if CMVSDE (\ref{SDE1}) does have  a weak solution $(X, Y)$ defined on some probability
space $(\O, \cF, \hP; \hF)$, then we can define a  new probability measure $\hQ$ so that 
$\frac{d\hQ}{d\hP}\big|_{\cF_t}=\bar L_t$, and $\bar L$ is a $\hP$-martingale. Moreover, if we define
the {\it unnormalized conditional distribution} by 
\bea
\label{unnorm}
\n_t(\cd):=\hE^\hQ[\bar L_t \1_{\{X_t\in \cd\}}|\cF^Y], 
\eea
then $\n$ should satisfy the nonlinear version of the {\it Zakai equation} (\ref{Zakai0}). In fact, one of the main purposes of this paper is to study the equivalence between the (weak) well-posedness of the SPDE (\ref{Zakai0})--(\ref{dual}) and that of the CMVSDE (\ref{SDE0}),
as was done in the case of McKean-Vlasov SDEs with common noise \cite{hammersley2019weak}.

In the section we first argue that the existence of the weak solution to the CMVSDE (\ref{SDE1}) implies the (weak) existence of the solution to SPDE (\ref{Zakai0})-(\ref{dual}). The converse, however, is a more difficult task, which will be one of the main results in this paper, to be presented in \S5. 

Let us begin by assuming that there exists a weak solution $(\O, \cF, \hP, \hF, B^1, B^2, X, Y)$  to the CMVSDE (\ref{SDE1}), and denote 
$\m=\hP^{X|Y}$. We define a process $\bar L$ such that
\bea
\label{barL}
\bar L_t=\exp\Big\{-\int_0^t h(s,X_s, Y_{\cd\wedge s}, \m_{\cd\wedge s})d B^2_s-\frac12\int_0^T|h(s,X_s, Y_{\cd\wedge s}, \m_{\cd\wedge s})|^2ds\Big\}, ~~t\in[0,T];
\eea
then $\bar L$ is an $\hF$-martingale under $\hP$, we can then define a new probability measure
$d\hQ=\bar L_Td\hP$, so that $(B^1, Y)$ is a $\hQ$-Brownian motion, and $(X, L)=(X, \bar L^{-1})$ is identical in law with the canonical solution to CMVSDE (\ref{SDE2}). Since we are interested only in 
conditional distribution, we can assume without loss of generality that $\O=\O^0$, $\cF=\cF^0$, and $\hQ=\hQ^0$. In what follows we simply
denote the weak solution by $(\hP, B^1, B^2, X, Y)$, and we call $(\hQ^0, 
X, L, \m=\hP^{X|Y}, B^1, Y)$ the ``canonical elements" of the given weak solution.
Now let us consider the unnormalized conditional distribution $\n_t(\cd)$, $t\in[0,T]$ defined by (\ref{unnorm}). For any bounded $ \f\in \hL^0(\hR)$, we denote
\bea
\label{condn}
\la \n_t, \f\ra:=\hE^{\hQ^0}[L_t \f(X_t)|\cF^Y_t], \qq t\in[0,T].
\eea
Furthermore, note that if $\psi:\hC_T\mapsto \hR$ is any bounded and measurable function, then
\beaa
\hE^{\hQ^0}[ L_t \f(X_t)\psi(Y_{\cd\wedge t})|\cF^Y_t]=\hE^{\hQ^0}[L_t \f(X_t)|\cF^Y_t]\psi(Y_{\cd\wedge t})=\hE^{\hQ^0}[ L_t \f(X_t)\psi(\by_{\cd
\wedge t})|\cF^Y_t]\big|_{\by=Y}.
\eeaa
A simple application of Monotone Class Theorem then leads to that, for any bounded measurable function $\f:[0,T]\times \hR\times\hC_T\mapsto \hR$, it
holds that
\bea
\label{eq1}
\hE^{\hQ^0}[ L_t \f(t,X_t, Y_{\cd\wedge t})|\cF^Y_t]&=&\hE^{\hQ^0}[ L_t \f(t, X_t, \by_{\cd
\wedge t})|\cF^Y_t]\big|_{\by=Y}\\
&=&\la \n_t, \f(t, \cd, \by_{\cd\wedge t})\ra|_{\by=Y}=\la \n_t, \f(t, \cd, Y_{\cd\wedge t})\ra. \nonumber
\eea
Our main result of this section is the following
\begin{prop}
\label{Zakaieq}
Assume Assumption \ref{assump1}. Let $(\hP, B^1, B^2, X,Y)$ be a weak solution to CMVSDE (\ref{SDE1}) and $(\hQ^0, X,  L, \m, B^1, Y)$ be its canonical elements. Let $\n_t$ be the unnormalized conditional distribution of $X_t$ given $\cF^Y_t$ under $\hQ^0$, defined by (\ref{condn}).
Then, for any $\f\in \hC^2_b(\hR^d)$, it holds that
\bea
\label{Zakai}
\la \n_t ,\f\ra&=&\la \n_0 ,\f\ra+\int_0^t\la \n_s ,\frac{1}{2}[\si\si^\top+\rho\rho^\top](s, \cd , Y_{\cd\wedge s},\mu_{\cd\wedge s}) :\nabla^2\f +b^\top(s, \cd , Y_{\cd\wedge s},\mu_{\cd\wedge s}) \nabla \f \ra ds\nonumber\\
&&+\int_0^t \la \n_s,(D\f, (\rho, h)^\top(s, \cd,Y_{\cd\wedge s}, \mu_{\cd\wedge s}))\ra dY_s, 
\eea 
where $A:B=\tr AB$, for $A, B\in\hS(\hR^d)$,  $D\f:=(\nabla \f, \f)^\top$ (see (\ref{Lcyl})), and $\m\in\hS_{\hF^Y}(\sP_1)$ is defined via
$\n$ by (\ref{dual}). 
\end{prop}

{\it Proof.} The proof is straightforward. Let $(\hQ^0, X,  L, \m, B^1, Y)$ be the canonical elements of the weak solution, then $(X, L)$ satisfies (\ref{SDE2}) under $\hQ^0$.  For any  {$\f\in \hC^2(\hR^d)$},  first applying It\^o's formula to $\f(X_t)$ and then to $L_t \f(X_t )$, 
then, denoting $\Th_t=( X_t, Y_{\cd\wedge t},\m_{\cd \wedge t})$, we have
\beaa
L_t \f(X_t )&=& \f(X_0)+\int_0^tL_s \{ \frac{1}{2}[\si\si^\top+\rho\rho^\top](s,   \Th_s):\nabla^2\f(X_s)+(\nabla \f)^\top(X_s)[\rho h](s,   \Th_s)\}ds\\
&&+\int_0^tL_s (\nabla \f(X_s))^\top dX_s+\int_0^tL_s \f(X_s)h(s,\Th_s)dY_s \\
&=& \f(X_0)+ \int_0^tL_s [\frac{1}{2}[\si\si^\top+\rho\rho^\top](s,   \Th_s):\nabla^2\f(X_s)+[ \tilde{b}+\rho h]^\top (s,   \Th_s)\nabla \f(X_s)]ds\\
&&+\int_0^tL_s( \si(s,  \Th_s)\nabla \f(X_s))^\top dB_s^1+\int_0^tL_s [\rho^\top (s, \Th_s)\nabla \f(X_s)+h(s,\Th_s)\f(X_s)]dY_s.
\eeaa 
Here we denote $\nabla \f$ and $\nabla^2 \f$ as the gradient and hessian of $\f$, respectively.
Taking conditional expectation $\hE^{\hQ^0}[\,\cd\,|\cF_t^Y]$ on both sides above, noting that $(B^1, Y)$ is a Brownian motion  under $\hQ^0$, it is readily seen (see, e.g., \cite{hammersley2019weak}) that, for any
$\hF^{B^1, Y}$-predictable process $H$, the following identities hold:
\bea
\label{Fubini}
\hE^{\hQ^0}\Big[\int_0^T H_t dB^1_t\ \Big|\cF^Y_t\Big]=0; \q 
\hE^{\hQ^0}\Big[\int_0^T H_t dY_t\ \Big|\cF^Y_t\Big]=\int_0^T\hE^{\hQ^0} [ H_t |\cF^Y_t]dY_t.
\eea
Since $\m=\hP^{X|Y}$ is $\hF^Y$-adapted, using the identity $ \tilde{b}+\rho h=b$, and
applying the two identities in (\ref{Fubini}) we obtain that  (suppressing variables)
\bea
\label{Zakai1}
&&\hE^{\hQ^0}[L_t f(X_t )|\cF_t^Y]\nonumber \\
\neg&\neg=\neg&\neg\hE^{ \hQ^0}[f(X_0 )|\cF_t^Y]
+ \int_0^t\hE^{\hQ^0}\big[L_s \big(\frac{1}{2}[\si\si^{\ts}+\rho\rho^{\ts}](\cd) :\nabla^2f(X_s)+b^{\ts}(\cd)\nabla f(X_s)\big)\big|\cF_s^Y\big]ds\nonumber\\
&&+\int_0^t\hE^{ \hQ^0}\big[L_s\big(\rho^{\ts}(\cd)\nabla f(X_s)+ h(\cd)f(X_s)\big)\big|\cF_s^Y\big]dY_s.
\eea
Now using   (\ref{eq1}) and the facts that $L_0=1$ and $\m$ is $\hF^Y$-adapted, we see that SPDE (\ref{Zakai}) is the same as 
(\ref{Zakai1}). 
\qed

\section{Superposition   for nonlinear Zakai equation and CMVSDE}
\label{relation 2}
\setcounter{equation}{0}

We now turn our attention to the converse of Proposition \ref{Zakaieq}, that is, the existence of the (weak) solution to the SPDE  (\ref{Zakai0})-(\ref{dual}) implies that of the CMVSDE (\ref{SDE1}). Compared to the usual CMVSDE with common noise case, the main difficulty here is the solution to the SPDE (\ref{Zakai0}) is the ``unnormalized" conditional distribution of $X$ given $\hF^Y:=\{\cF^Y_t\}$ under the reference probability $\hQ^0$, thus a ``total probability" type conditioning argument does not lead to the joint $\hQ^0$-distribution of $(X, Y)$, despite the fact that $Y$ is a $\hQ^0$-Brownian motion. However, the idea of the so-called ``superposition" is still effective. 

To facilitate our argument, in what follows we strengthen our standing assumption (Assumption \ref{assump1}) slightly by adding the following condition on the coefficients $b, h, \rho$. 
\begin{assum}
\label{assump2}
For fixed $(t, \by, \m)\in [0,T]\times\hC_T\times \sC_T(\sP_1)$, the function $\f(t, \cd, \by, \m)\in  \hC^{(1)}_b(\hR^d)$, where $\f=b,  h, \rho$. 
In particular, we shall assume that the $\hC^{(1)}$ norm of the functions $b, h, \rho$ is the same as the Lipschitz constant $C>0$ in 
Assumption \ref{assump1}, uniformly in $(t, \by, \m)$.
\qed
\end{assum}

 {Let us recall the canonical space $(\O^0, \cF^0, \hF^0, \hQ^0)$ and the (canonical) Brownian motion $(B^1, Y)$. Recall also the space $\hS^2_{\hF^Y}(\sP_1)$ defined by (\ref{SpG}), where
$\sP_1=\sP_1(\O^0)$.
We have the following result.
\begin{thm}
 \label{superpos1}
 Assume that Assumptions \ref{assump1} and \ref{assump2} are in force. Assume that there exists a measure-valued process $\n\in \hS^2_{\hF^Y}(\sM_1)$,   defined on the canonical space, such that the SPDE (\ref{dual})--(\ref{Zakai0}) hold. Then, there exists a probability measure $\hP^\n\in \sP_1(\O^0)$, and  $\hF^0$-progressively measurable processes $(X^\n, L^\n, B^{2, \n})\in\hL^2_{\hF^0}(\O;\hC([0,T];\hR^{d+2}))$, such that,

 (i) $(X^\n, L^\n)$ satisfies SDE (\ref{SDE2}), $\hP^\n$-a.s.; 
 
 (ii)   $\frac{d\hP^\n}{d\hQ^0}\big|_{{\cF}_T}=L^\n_T$;
 
 (iii)  $(B^1, B^{2,\n})$ is a $(\hP^\n, \hF^0)$-Brownian motion;
 
 (iv) $(\O^0, \cF^0, \hP^\n, \hF^0, B^1, B^{2,\n}, X^\n, Y)$ is a weak solution to CMVSDE (\ref{SDE1}). 

 \end{thm}
 }
 
 {\it Proof.} 
 Throughout this proof we assume that $\n$ is a solution for the nonlinear Zakai equation defined in \eqref{Zakai0}, and we try to construct a weak solution  $(\tilde{\O}, \tilde{\cF}, \hP, \tilde{\hF}, B^1, B^2, X, Y)$ of   SDE (\ref{SDE2}), such that the conclusions (i)-(iv) hold. 

We begin by considering the canonical  space $(\O^0, \cF^0, \hF^0, \hQ^0)$ and let $\n\in \hS^2_{\hF^Y}(\sM_1)$ be the solution of the SPDE ((\ref{dual})--\ref{Zakai0})). Let $(B^1, Y)$ be the canonical Brownian motion on $(\O^0, \cF^0,\hQ^0)$, and consider the following SDE  with random coefficients:
\bea
\label{SDE-random-0}
\left\{\ba{lll}
dX^{\n}_t=\tilde{b}(t,  X^{\n}_t, \o^2_{\cd\wedge t},  {\m}^\n_{\cd\wedge t})dt+ \si(t,  X^{\n}_t, \o^2_{\cd\wedge t},  \m^\n_{\cd\wedge t})dB^{1}_t+\rho(t, X^{\n}_t, \o^2_{\cd\wedge t},  \m^\n_{\cd\wedge t})dY_t;\ms\\
dL_t^{\n}= h(t,  X^{\n}_t, \o^2_{\cd\wedge t},  \m^\n_{\cd\wedge t})L^\n_t dY_t, \q\qq  X_0=x, \q L_0=1.
\ea\right.
\eea 
Here in the above, $\m=\m^{\n}\in \hS^2_{\hF^Y}(\sP_1)$ is defined by (\ref{dual}),  $\tilde b:=b-\rho h$, and for $\phi =b,\si, \rho, h$, by a slight abuse of notations we identify
$$ \phi(t, x, \o, \m):=\phi(t, x, \o^2_{\cd\wedge t}, \m), \qq t\in[0, T], \q\o=(\o^1, \o^2)\in \hC^d_T\times\hC_T, \m\in\sC_T(\sP_1),
$$
so that $\tilde b, \si, \rho$ and $h$ are all $\hF$-progressively measurable. Consequently, as now an SDE with random coefficients \reff{SDE-random-0} has a unique (strong) solution on $(\O^0, \cF^0, \hQ^0)$, thanks to Assumption \ref{assump1}. Furthermore, since $h$ is bounded, we know that $L_T^{\n}$ is a $\hQ^0$-martingale, and $d\hP:= L^{\n}_Td\hQ^0$ defines a probability measure such that $(B^1,B^{2,\n})$ is an $(\hF,\hP)$-Brownian motion, where $B^{2,\n}_t:=Y_t-\int_0^th(\cds)ds$, $t\in[0,T]$. Now note that the under $\hQ^0$, the  canonical process $\o^2_{\cd\wedge t}=Y_{\cd\wedge t}$, $t\in[0,T]$. We shall now replace $\o^2$ in the SDE \reff{SDE-random-0} by $Y$. Then, under 
$\hP$, the $(X^\n, Y, B^1, B^{2,\n})$ should satisfy the SDE:
\bea
\label{SDE-random-1} 
\left\{\ba{lll}
dX_t^{\n}=b(t,  X_t^{\n}, Y_{\cd\wedge t},  {\m}_{\cd\wedge t})dt+ \si(t,  X_t^{\n}, Y_{\cd\wedge t},  {\m}_{\cd\wedge t})dB^1_t+\rho(t, X_t^{\n}, Y_{\cd\wedge t},  {\m}_{\cd\wedge t})dB^{2,\n}_t;\ms\\
dY_t = h(t,  X_t^{\n}, Y_{\cd\wedge t},  {\m}_{\cd\wedge t})dt+dB^{2,\n}_t,  \qq \q  X_0=x, \q Y_0=0.
\ea\right.
\eea 

Now, note that $\m=\m^{\n}\in \hS_{\hF^Y}^2(\sP_1)$, by Doob-Dynkin theorem, there exists an $\hF^Y$-progressively measurable Borel functional $\Phi^\n:\hC_T \mapsto \scC(\mathcal P_1)$, such that  
\bea
\label{Phi map}
\m_t=\m_t^\n =\Phi^\n(Y)_t=\Phi^\n(Y_{\cdot\wedge t})_t,\quad t\in [0,T],\quad \hP\as
\eea
Next, let us define: for $\phi=b, \si, \rho, h$, 
\bea 
\label{Phinu} 
\phi^\n(t,  x, \by )=\phi\circ \Phi^\n(t,  x, \by):=\phi(t,x, \by, \Phi^\n(\by)_t), \q (t, x, \by)\in[0,T]\times \hR\times \hC_T, 
\eea 
Then, the SDE (\ref{SDE-random-1}) can be written as
\bea
\label{SDEnu}
\left\{\ba{lll}
dX^{\n}_t=b^\n(t,  X^{\n}_t, Y_{\cd\wedge t})dt+ \si^\n(t,  X^{\n}_t, Y_{\cd\wedge t})dB^{1}_t+\rho^\n(t, X^{\n}_t, Y_{\cd\wedge t})dB^{2}_t;\ms\\
dY_t^{\n}= h^\n(t,  X^{\n}_t, Y_{\cd\wedge t})dt + dB^{2,\n}_t, \q\qq  X_0=x, \q Y_0=0,
\ea\right.
\eea 
and $(X^\n, Y^\n)$ is a (pathwisely) unique strong solution for a given $\n\in\hS^2_{\hF^Y}(\sM_1)$. We now consider SDE \reff{SDEnu},  as a standard filtering problem. We first note that by 
definition of $\hP^\n$ we can write $d\hQ^0=\hat L^\n_Td\hP^\n$, where $\bar L^\n:=(L^\n)^{-1}$ satisfies the SDE
\beaa
d \bar L^\n_t= -h^\n(t,   X^\n_t, Y^\n_{\cd\wedge t}) \bar L^\n_tdB^{2,\n}_t,  \qq \q  L^\n_0=1, \q t\in[0,T],
\eeaa
and $(B^1, Y^\n)=(B^1, Y)$ is a $\hQ^0$-Brownian motion. Let us now define, for $A\in\sB(\hR^d)$, 
 \bea 
 \label{conditional mu bar}
 \bar \m^2_t(A):=(\hP^\n)^{X^\n|Y}(A):=\hP^\n[ X_t^{\n}\in A | \cF_t^{Y} ]=\frac{\hE^{\hQ^0}[\bar L_t^\n\1_{\{X_t^{\n}\in A\}}|\cF_t^{Y}]}{\hE^{\hQ^0}[\bar L_t^{\n}|\cF_t^{Y}]}:=\frac{\bar\n_t( A)}{\bar\n_t(\hR^d)}, 
 \eea 
where $\bar\n_t(\cd):= \hE^{\hQ^0}[\bar L_t^\n\1_{\{X_t^{\n}\in \cd\}}|\cF_t^{Y}]$, $t\in[0,T]$, is the corresponding unnormalized conditional distribution. Thus, by the standard nonlinear filtering theory, $\bar\n(\cd)=(\bar\n_t(\cd))_{t\in[0,T]}$ satisfies the (linear) Zakai equation under $\hQ^0$: for any $\f\in\hC^2(\hR^d)$, 
\bea
\label{Zakai-l}
\la \bar \n_t ,\f\ra&=&\la \bar\n_0 ,\f\ra+\int_0^t\la \bar\n_s ,\frac{1}{2}[\si^\n(\si^\n)^{\top}+\rho^\n(\rho^\n)^{\top}](s, \cd ,\o^2_{\cd\wedge s}) :\nabla^2\f +(b^\n)^{\top}(s, \cd , \o^2_{\cd\wedge s}) \nabla \f \ra ds\nonumber\\
&&+\int_0^t \la \bar\n^2_s,(D\f, (\rho^\n, h^\n)^\ts(s, \cd, \omega^2_{\cd\wedge s}))\ra dY_s.  
\eea
Note that the equation (\ref{Zakai-l}) is a linear, measure-valued SPDE with random coefficients (by identifying the canonical process $Y_t(\o)=\o^2_t$), where coefficients $b^\n, \si^\n, \rho^\n$ and $h^\n$ are defined by \reff{Phinu}. Since $\n$ is a solution to \reff{Zakai0}, by definition it is also a solution to (\ref{Zakai-l}).
 It is now suffices to show the uniqueness of the linear equation \eqref{Zakai-l}, which then implies 
$\n^2=\bar\n^2$, and consequently, the solution $(X^\n, Y^\n)$ is a weak solution to CMVSDE under $\hP^\n$.  

The proof of the uniqueness of (\ref{Zakai-l}) follows from a similar argument in \cite{kurtzxiong, lacker2020}. We give a brief sketch for completeness.
To begin with, for $\delta>0$, denote the Gaussian kernel $G_{\delta}(x):=(2\pi\delta)^{-d/2}e^{-|x|^2/(2\delta)}$, and for   $\n\in\hS^2_{\hF^Y}(\sM_1)$, we define a random field $T^{\n,\d}:[0,T]\times\O\times\hR^d\to\hR$ by
\bea
\label{Zdn}
T^{\n,\d}_t (x):=(G_{\delta}*\n_t)(x)=\int_{\hR^d}G_{\delta}(x-y)\n_t(dy),\q   (t,x)\in [0,T]\times\hR^d, ~\hP\as
\eea
Clearly,  $T^{\n,\delta}_t\in \hL^2(\hR^d)$, $t\in[0,T]$, $\hP$-a.s., and for any $\f\in\hC^{2}(\hR^d)$, and
 $\n$ satisfying the linear  Zakai equation (\ref{Zakai-l}), we can easily check, using Fubini theorem, integration by parts, as well as the fact that $G_\delta$ and differentiation commute,  that
\beaa 
\la T^{\n,\d}_t,\f\ra_2  &= &\la T^{\n,\d}_0,\f\ra_2+\int_0^t\big\la\frac{1}{2}\sum_{i,j=1}^d\pa_{ij}(G_{\delta}* [( \si\si^{\top}\neg+\neg\rho\rho^{\top})_{ij}\n_s ] )-\sum_{i=1}^d\pa_i(G_{\delta}*(b^i\n_s) ),\f \big\ra_2 ds\\
&&+\int_0^t\la G_{\delta}*(h_s\n_s)-\sum_{i=1}^d\pa_i(G_{\delta}*\rho_s^{i} \n_s) ,\f \ra_2    dY_s,
\eeaa
where $\la\cd,\cd\ra_2$ denotes the inner product of $\hL^2(\hR^d)$. Next,  we define, for $t\in[0,T]$,
\bea
 \label{alpha}
 \alpha_t= \|b(t, \cdot)\|_{\hL_\infty(\O\times C^1_b)}+\|h(t,\cd )\|_{\hL_\infty(\O\times C_b^1)}+\|\rho(t, \cdot)\|_{\hL_\infty(\O\times C_b^2 )}^{2}+\|\sigma(t, \cdot)\|_{\hL^\infty(\O\times C_b^2)}^{2}. 
 \eea 
 Then, $\a \in\hL_\infty([0,T])$, thanks to Assumption \ref{assump2}. Furthermore, define
 $\beta_t:=e^{-K\int_0^t \alpha_sds}$, where the constant $K$ is to be determined. and applying   
It\^o's formula we have
 \bea
 \label{betaTn}
d[\beta_t\la T_t^{\n,\delta},\f \ra^2_2]\neg&\neg\neg=\neg\neg&\neg\beta_t\big[-  K \alpha_t\la T_t^{\n,\delta},\f\ra^2_2+\big\la G_{\delta}*(h_t\n_t)-\sum_{i=1}^d\pa_i(G_{\delta}*\rho_s^{i} \n_t) ,\f \big\ra^2_2 \nonumber\\
 &&+ 2 \la T_t^{\n,\delta},\f\ra_2 \big\la\frac{1}{2}\sum_{i,j=1}^d\pa_{ij}(G_{\delta}* [( \si\si^{\top}+\rho\rho^{\top})^{ij}\n_t ] )-\sum_{i=1}^d\pa_i(G_{\delta}*(b^i\n_t) ),\f \big\ra_2 \big]dt\nonumber\\
 &&+  2\beta_t \la T_t^{\n,\delta},\f\ra \big\la G_{\delta}*(h_t\n_t)-\sum_{i=1}^d\pa_i(G_{\delta}*\rho_t^{i} \n_t) ,\f \big\ra_2    dY_s.
 \eea 
Let $\{\f^i\}_{i=1}^\infty$ be an orthornormal basis of $\hL^2(\hR^d)$. Applying \reff{betaTn} for each $\f^i$ and summing up, then taking expectation under $\hQ^0$, we obtain from 
Parseval's identity that
 \bea \label{ito1}
&&d\big( \beta_t\hE^{\hQ^0}\big[\| T_t^{\n,\delta}\|_2^2\big]\big)=\beta_t\big\{\neg-\neg K \alpha_t\hE^{\hQ^0}\big[\| T_t^{\n,\delta}\|_2^2\big]+\hE^{\hQ^0}\big[\| G_{\delta}*(h_t\n_t)-\sum_{i=1}^d\pa_i(G_{\delta}*\rho_t^{i} \n_t) \|_2^2\big]\big\}dt\nonumber \\
 &&\qq+ \beta_t \hE^{\hQ^0}\big[\big\la\sum_{i,j=1}^d\big[\pa_{ij}(G_{\delta}* [( \si\si^{\top})^{ij}
 +(\rho\rho^{\top})^{ij}\big]\n_t  ),T_t^{\n,\delta} \big\ra \big]dt \\
 &&\qq - 2\beta_t\hE^{\hQ^0}\big[ \big\la\sum_{i=1}^d\pa_i(G_{\delta}*(b^i\n_t) ),T_t^{\n,\delta} \big\ra \big]dt.
\nonumber
  \eea 
Now, note that by definition $\|T^{\n, \d}_t\|_2 \le \|G_\d*|\n_t|\|_2$,  where $|\n|$ denotes the total variation measure of $\n$. Therefore, applying \cite[Lemma 3.2 \& 3.3]{kurtzxiong}, we have
$$\left\{\ba{lll}
\dis | \la T_t^{\n,\delta},\sum_{i=1}^d\pa_i(G_{\delta}*(b^i\n_s) ) \ra|  \le  C \|b(t, \cdot)\|_{\hL_\infty(\O\times C_b^1)}\|G_{\delta}*|\n_t| \|_2^2 
\ms\\
{\small \dis |\big\la\sum_{i,j=1}^d\neg\neg \big[\pa_{ij}(G_{\delta}* [( \si\si^{\ts})^{ij}
\neg+\neg(\rho\rho^{\ts})^{ij}\big]\n_t  ),T_t^{\n,\delta} \big\ra|\le  C \big(\|\si(t, \cdot)\|^2_{\hL_\infty(\O\times C_b^2) }\neg+\neg\|\rho(t, \cdot)\|^2_{\hL_\infty(\O\times C_b^2) })\|G_{\delta}*|\n_t| \|_2^2,}
  \ea\right.
$$
where $C>0$ is a generic constant that may depend on $d$ but independent of the coefficients, and is allowed to vary from line to line. With the similar estimate on the remaining term on the right hand side of  \eqref{ito1}  we deduce  the
following estimate: 
 \bea 
 \label{decreasing condition}
 \frac{d}{dt}\big( \beta_t \hE^{\hQ^0}\big[ \| T_t^{\n,\delta} \|_2^2\big]\big)&\le& C\beta_t\a_t \hE^{\hQ^0}\big[\|G_{\delta}*|\n_t| \|_2^2  \big]-K\beta_t \alpha_t \hE^{\hQ^0}\big[\| T_t^{\n,\delta}\|_2^2\big]\\ 
 &\le&
 (C-K)\beta_t\alpha_t\hE^{\hQ^0}\big[ \|G_{\delta}*|\n_t| \|_2^2\big]+K\beta_t \alpha_t \hE^{\hQ^0}\big[ \|G_{\delta}*|\n_t| \|_2^2-\| T_t^{\n,\delta}\|_2^2\big]\nonumber\\
 &\le&K\beta_t \alpha_t \hE^{\hQ^0}\big[ \|G_{\delta}*|\n_t| \|_2^2-\| T_t^{\n,\delta}\|_2^2\big].\nonumber
 \eea 
In the above we choose the constant $K\ge C$ in (\ref{alpha}). We note that if $\bar\n^1$ and $\bar\n^2$ are two (positive) measure-valued processes satisfying the linear  Zakai equation (\ref{Zakai-l}), then we have $T^{\n^i,\delta}_t(x):=(G_{\delta}*\bar \n^i_t)(x)$, $i=1,2$, so that $\frac{d}{dt}\hE^{\hQ^0}\Big[  \beta_t{\| T^{\n^i,\delta}_t\|_2^2}\Big]\le 0$, for $i=1,2$, thanks to \eqref{decreasing condition}, or $\hE^{\hQ^0}\big[  \beta_t\| T^{\n^i,\delta}_t\|_2^2 \big]\le \hE^{\hQ^0}\big[\| T^{\n^i,\delta}_0\|_2^2\big]\le \hE^{\hQ^0}\big[\| \rho^i\|_2^2\big]$. We can then follow  a line-by-line argument  as that in \cite[Proposition 3.2]{lacker2020}, to first argue that both $\n^1$ and $\n^2$ have densities $\rho^i=\{\rho^i_t\}\in \hL^2([0,T]\times\hR^d)$, and then applying the dominated convergence theorem to conclude that $\frac{d}{dt}\hE^{\hQ^0}\Big[  \beta_t{\| \rho^1_t-\rho^2_t\|_2^2}\Big]\le 0$, which implies that the linear 
 Zakai equation with random coefficient \eqref{Zakai-l} has a unique solution $\n$, which has a density $\rho=\{\rho_t\}\in\hL^2([0, T]\times\hR^d)$. This complete the proof of uniqueness of the linear Zakai equation (\ref{Zakai-l}), whence the theorem.
\qed

\section{Superposition for CFPE and nonlinear Zakai equation}
\label{relation 1}
\setcounter{equation}{0}

In this section we turn our attention to our second superposition principle, regarding the non-linear Zakai equation  (\ref{Zakai}). More precisely, we shall consider (\ref{Zakai}) as an SDE taking values in the space 
 $ \sM_1(\hR^d) $ and consider its corresponding Fokker-Planck equation (CFPE, see \eqref{FK0}), defined on the probability laws of $ \sM_1(\hR^d) $-valued processes, which we shall refer to as the {\it conditional Fokker-Planck equation}. We shall argue that a solution to   the conditional Fokker-Planck equation can be ``lifted" to a weak solution to the SPDE (\ref{Zakai}) by the superposition principle. For technical reasons, starting from this part  we shall assume that the dependence on the (normalized) conditional law $\m$ is ``state dependent" instead of ``path-dependent", that is, depending on $\m_t$ instead of $\m_{\cdot\wedge t}$ in all the coefficients. This simplification allows us to relate the solution to the conditional Fokker-Planck equation for the non-linear Zakai equation via the simple Bayes' rule. 

To be more specific, let us first recall the 
non-linear Zakai equation \eqref{Zakai} with state dependence on $\m$, defined on the canonical space $(\O^0, \cF^0, \hQ^0)$: for any 
$\f\in\hC^2_b(\hR^d)$, 
\bea
\label{Zakai: state}
\la \n_t ,\f\ra&=&\la \n_0 ,\f\ra+\int_0^t\la \n_s ,\frac{1}{2}[\si\si^{\top}+\rho\rho^{\top}](s, \cd , Y_{\cd\wedge s},\mu_{s}) :\nabla^2\f +b^{\top}(s, \cd , Y_{\cd\wedge s},\mu_{ s}) \nabla \f \ra ds\nonumber\\
&&+\int_0^t \la \n_s,(D\f, (\rho, h)^\top (s, \cd,Y_{\cd\wedge s}, \mu_{s}))\ra dY_s, 
\eea 
where $Y$ is the canonical process, whence a $\hQ^0$-Brownian motion,  and 
$\n= \{\n_t\}$ is the unnormalized conditional probabilities defined in the previous sections, so that  $\n_t(\o)\in \sM(\hR^d)$, $\m_t(\o)\in\sP_1(\hR^d)$, $t\in[0,T]$,  $\hQ^0$-a.e. $\o\in\O^0$, and that $\la \m_t, \f\ra=\frac{\la v_t, \f\ra}{\la v_t, 1\ra}$.
  
Now for $t\in[0,T]$, let us recall \eqref{rcpd-y} the  {\it regular conditional probability distribution} (RCPD) of $\hQ^0$ given $\cF^Y_t$, on $(\O^0, \cF^0)$, denoted by $\hP^y_t(\cd):=\hE^{\hQ^0}[\, \cd \, |\cF^Y_t](y)\in \sP(\O^0)$, for $\hQ^0$-a.e. $y\in\O^0$. 
Consider  \eqref{Zakai: state} as an SDE for $\la\n_t, \f\ra $ on $\sM(\hR^d)$ with $\f\in\hC^2_b(\hR^d)$. For each $t\in[0,T]$, let us define
the conditional law of $\n_t$ under $\hP^Y_t$ by
\bea
\label{bPyt}
\bP_t^Y:=\hP^Y_t\circ (\n_t)^{-1}\in \sP_1(\sM_1(\hR^d)), \qq \hQ^0\textit{-a.s.}, 
\eea
and define a mapping $\bm:\sM_1(\hR^d)\mapsto \sP_1(\hR^d)$, by
\bea
\label{bm}
\la \bm(\bn),\f\ra =\frac{\la \bn,\f\ra}{\la \bn,1\ra}, \qq \bn\in\sM_1(\hR^d).
\eea
We first show that for $\hQ^0$-a.s., $\bP_t^Y$, $t\in[0,T]$,  satisfies the conditional Fokker-Planck equation \eqref{FK0} or equivalently \eqref{FP L}. 
Then, we shall prove the superposition principle between the conditional Fokker-Planck equation and the nonlinear Zakai equation (\ref{Zakai}). 
We should note that the superposition principle regarding the Zakai-type SPDE has been studied for the ``common noise" case (i.e., $h\equiv0$) in recent  works \cite{barbu2020nonlinear, lacker2020, Qiao2022, ren2020space}. However, to the best of our knowledge, the superposition principle involving the conditional
Fokker-Planck equation (i.e. stochastic version of Fokker-Planck equation) corresponding to the general form such as 
\eqref{Zakai: state} is new. 

To simplify the argument, in what follows we shall make use of the following  assumption.
\begin{assum}
\label{assump2}
The coefficients $\f=(b,\si, \rho, h):[0,T]\times \hR^d\times \hC_T\times \sM_1\mapsto \hR^d\times \hR^{d\times d}\times \hR^{d\times d}\times \hR$ are bounded and countinuous, such that 
for fixed $(t, \by, \m)\in [0,T]\times \hC_T\times \sM_1$, the mapping $x\mapsto \f(t, x, \by, \m)$ is twice continuously differentiable with bounded derivatives, uniformly in $(t, \by, \m)$. 
\qed
\end{assum}

Let us begin with the following straightforward result.
\begin{prop}\label{prop: fokker plack}
Assume that Assumption \ref{assump2} is in force. Let $(\Omega,\cF,\hQ^0)$ be a probability space on which is defined a $d+1$-dimensional Brownian motion $(B^1,Y)$ and a continuous $\hF$-adapted $\sM_1(\hR^d)$ valued process $(\n_t)_{t\in[0,T]}$ satisfying SPDE \eqref{Zakai: state}. 
Then, the $\sP_1(\sM_1(\hR^d))$-valued process $\{\bP^Y_t\}_{t\in[0,T]}$ defined by (\ref{bPyt}) satisfies the following conditional Fokker-Planck equation: 
\bea
\label{CFPeq}
\int_{\sM(\hR^d)}F (\bn)(\bP^Y_t-\bP^Y_0)(d\bn)
&=&\int_0^t\int_{\sM(\hR^d)} \cL_s^{Y, \bm} F(\bn) \bP^Y_s(d\bn) ds \\
&&+\int_0^t\int_{\sM(\hR^d)}\bd_{\bn}F(\bn,x)\cdot \cH(t, x , Y_{\cdot\wedge s},{\bm}) \bP^Y_s(d\bn) dY_s, \qq \hQ^0\textit{-a.s.}\nonumber
\eea
where the operator $\cL^{Y, \bm}_t$ is defined as  in \eqref{FP operator}, and $F$ is a $k$-cylindrical test function. 
\end{prop}

{\it Proof.} Let $\n=(\n_t)$ be the $\hF$-adapted, $\sM_1(\hR^d)$ valued continuous process that satisfies \eqref{Zakai: state}.
Let $k\in\hN$, and $\psi_1, \cds, \psi_k\in\hC^\infty_c(\hR^d)$, and $f\in\hC^2(\hR^k)$.  Consider  the $k$-{\it cylindrical} test function $F\in\hC^2_b(\sM_1)$ defined by 
\bea 
\label{Fn}
F(\bn)=f(\la \bn, \psi\ra)=f(\la\bn,\psi_1\ra, \cdots, \la\bn,\psi_k\ra ), \qq \bn\in\sM_1(\hR^d).
\eea 

Applying It\^o's formula to  $f(\la \n_t, \psi\ra)=f(\la \nu_t, \psi_1\ra,\cdots,\la \nu_t, \psi_k\ra)$ we have
\bea
\label{ito expansion}
df(\la \nu_t, \psi\ra)\neg\neg&\neg\neg=\neg\neg&\neg\neg\sum_{i=1}^k\pa_i f(\la \nu_t, \psi\ra)\Big[\la \n_t,(\nabla \psi_i(\cd))^{\top}\rho(t, \cd,Y_{\cd\wedge t}, \mu_{t})+h(t, \cd, Y_{\cd\wedge t}, \mu_{t})\psi_i(\cd)\ra\Big] dY_t\nonumber\\
&&\neg+ \sum_{i=1}^k\pa_i f(\la \nu_t, \psi\ra) \Big[\la \n_t ,\frac{1}{2}\si\si^{\top}(t, \cd , Y_{\cd\wedge t},\mu_{t}):\nabla^2 \psi_i(\cd)+(\nabla\psi_i )^{\top}b(t, \cd , Y_{\cd\wedge t},\mu_{t})\ra \Big]dt\nonumber\\
&&\neg+\frac{1}{2}\sum_{i,j=1}^k\pa_{ij} f(\la \nu_t, \psi\ra)\Big[\la \n_t,(\nabla \psi_i(\cd))^{\top}\rho(t, \cd,Y_{\cd\wedge t}, \mu_{t}) +h(t, \cd, Y_{\cd\wedge t}, \mu_{t})\psi_i(\cd)\ra\Big] \\
&&\qq\qq\times \Big[\la \n_t,(\nabla \psi_j(\cd))^{\top}\rho(t, \cd,Y_{\cd\wedge t}, \mu_{t})+h(t, \cd, Y_{\cd\wedge t}, \mu_{t})\psi_j(\cd)\ra\Big] dt.\nonumber
\eea
Now recall that $\bP^y_t:=\hP^y_t\circ (\n_t)^{-1}\in \sP_1(\sM_1(\hR^d))$ and $\hP^y_t(\cd):=\hE^{\hQ^0}[\, \cd \, |\cF^Y_t](y)\in \sP(\O^0)$ is the {\it regular conditional probability distribution} (RCPD) on $(\O^0, \cF^0)$, and we define $\bP^Y_t(y):=\bP^y_t$, for $(t,y)\in [0,T]\times \Omega_0$. Integrating \eqref{ito expansion} on time interval $[0,t]$, taking conditional expectation $\hE[\cdot |\cF_t^Y ]$ under measure $\hQ^0$ (under which $Y$ is a Brownian motion), and applying Fubini theorem for conditional expectations, for $\f \in \hL^2_{\hF}([0,T])$, we have  $\hE[\int_0^t \f(s) ds |\cF_t^Y]=\int_0^t \hE[\f(s) |\cF_s^Y]ds$, and 
\beaa 
\hE\left[\int_0^t \f(s) dY_s \Big|\cF_t^Y\right]=\int_0^t \hE[\f(s) |\cF_t^Y]dY_s=\int_0^t \hE[\f(s) |\cF_s^Y\vee \cF^Y_{s,t}]dY_s=\int_0^t \hE[\f(s) |\cF_s^Y]dY_s,
\eeaa 
where the last equality above is due to the fact that $\cF^Y_s\ind \cF^Y_{s,t}$ since $Y$ has independent increments.
We obtain the following conditional Fokker-Planck equation for $\bP^Y_t$,    for $t\in[0,T]$   and $\hQ^0$ a.s.,
{\small
\bea
\label{FP for nonlinear Zakai}
&&\int_{\sM(\hR^d)}f(\la \bn,\psi\ra)(\bP^Y_t-\bP^Y_0)(d\bn)\nonumber \\
&=&
\int_0^t\int_{\sM(\hR^d)}\Big\{\sum_{i=1}^k\pa_if(\la \bn, \psi\ra) \Big[\la \bn ,\frac{1}{2}\si\si^{\top}(s, \cd , Y_{\cdot\wedge s},{\bm}):\nabla^2 \psi_i(\cd)+(\nabla \psi_i(\cd))^{\top}b(s, \cd , Y_{\cdot\wedge s},{\bm}) \ra \Big]\Big\}\bP^Y_s(d\bn)ds \nonumber\\
&&+\int_0^t\int_{\sM(\hR^d)}\Big\{\frac{1}{2}\sum_{i,j=1}^k\pa_{ij}f(\la \bn, \psi\ra)\Big[\la \bn,(\nabla\psi_i)^{\top}\rho(s, \cd,Y_{\cdot\wedge s}, {\bm})+h(s, \cd, Y_{\cdot\wedge s},{\bm})\psi_i(\cd)\ra\Big] \\
&&\qq\qq\qq\qq\qq\qq\qq\times \Big[\la \bn,(\nabla\psi_j)^{\top}\rho(s, \cd,Y_{\cdot\wedge s}, {\bm})+h(s, \cd, Y_{\cdot\wedge s},{ \bm})\psi_j(\cd)\ra\Big]
\Big\}\bP^Y_s(d\bn) ds,\nonumber\\
&&+\int_0^t\int_{\sM(\hR^d)} \Big\{\sum_{i=1}^k\pa_i f(\la \bn, \psi\ra)\Big[\la \bn,(\nabla \psi_i(\cd))^{\top}\rho(t, \cd,Y_{\cdot\wedge s}, \bm)+h(t, \cd, Y_{\cdot\wedge s}, \bm)\psi_i(\cd)\ra\Big]  \Big\}\bP^Y_s(d\bn)dY_s.\nonumber
\eea
}
Rewriting the above expressions by using the $L$-derivatives (Definition \ref{L derivative}) and noting the definitions (\ref{Lcyl}), (\ref{Fn}) and 
that $\bm=\bm(\bn)$, we derive (\ref{CFPeq}).  
\qed 
\begin{rem}
Note that, if $h=0$, we could define
\beaa 
\bd_{\bn}^2F(\bn,x,x')=\sum_{i,j=1}^k\pa_{ij}f(\la\bn,\psi\ra)\nabla\psi_i(x)\nabla\psi_j(x')^{\top},
\eeaa 
which produces the Fokker-Plance equation on $\sP(\hR^d)$ in \cite{lacker2020}. If $h\neq 0$, but the coefficients depend only on the signal process $X$,  the corresponding Fokker-Planck equation is derived in \cite{Qiao2022} for Zakai equation with common noise, which is not of the Mckean-Vlasov type as in the present case. 
\qed
\end{rem}

 Our main goal of this section is to prove the superposition for the FPE (\ref{CFPeq}) and the nonlinear Zakai equation \reff{Zakai: state}, that is, 
  the existence of the solution to the conditional Fokker-Planck equation (\ref{CFPeq}) would imply the existence of the (weak)  solution to \reff{Zakai: state}. 
It is  worth noting that a similar result has been studied in \cite[Theorem 1.5]{lacker2020} when only the common noise is considered, i.e. $h=0$. Furthermore, the existence of solution for the Fokker-Planck equation with $h\equiv0$ has been studied in \cite{barbu2020nonlinear}.  
However, to the best of our knowledge, the case with $h\neq 0$ has yet been studied in the literature.

 For technical reasons in what follows we shall consider only the Markovian case on the 
 observation process $Y$. That is we shall reduce the path-dependence of the coefficients (on $Y_{\cdot\wedge t}$ to the state-dependent (on $Y_t$) so that the set of coefficients $\f=(b,\si, \rho, h)$ in Assumption \ref{assump2} becomes the mapping from $[0,T]\times \hR^d\times \hR\times \sM_1$ to $\hR^d\times \hR^{d\times d}\times \hR^{d\times d}\times \hR$ . 
 This simplification is mainly to facilitate our argument using the so-called {\it filtered martingale problem} in proving the superposition principle on $\mathbb R^{\infty}$ associated with the non-linear Zakai equation. We hope to be able to address the general path-dependent case in our  future studies. 
Let us first specify the precise meaning of the solution to the conditional Fokker-Planck equation (\ref{FP for nonlinear Zakai}). 
 \begin{defn}
 \label{FPsol}
A family of $\sP_1(\sM_1(\hR^d)))$-valued process $(\bP_t^Y)_{t\in[0,T]}$  is called a solution to the conditional Fokker-Planck equation (\ref{FP for nonlinear Zakai}), if, 

\ms
(i) $(\bP_t^Y)_{t\in[0,T]}\in \hL^2_{\hF^Y}(\Omega_0;  \hC([0,T],\sP_1(\sM_1(\hR^d))))$; and 

\ms
(ii) for any  $k$-cylindrical test function $F(\cd)$ defined by (\ref{Fn}), the equation (\ref{CFPeq}), or equivalently (\ref{FP for nonlinear Zakai}), holds. 

In particular, the solution $(\bP_t^Y)$ is called ``regular" if for some $p>1$ and $\hQ^0$-a.s., the following integrability condition holds: 
\bea
\label{integ}
&&\Phi(T):=\int_0^T \int_{\sM(\hR^d)}\Big[ \|b(t,\cd, Y_t ,\bm)\|_{L^1(|\bn|)}^p+\|\sigma\sigma^{\ts}(t,\cd, Y_t,\bm)\|_{L^1(|\bn|)}^p\\
&&\qq\qq\qq+\|\rho(t,\cd, Y_t,\bm)\|_{L^1(|\bn|)}^{2p}+\|h(t,\cd, Y_t,\bm)\|_{L^1(|\bn|)}^{2p}\Big] \bP^Y_t(d\bn)dt<\infty, \nonumber
\eea  
where $\bm=\bm(\bn)$ is defined by (\ref{bm}), and $\Phi(\cdot)$ is $\hF^Y$-adapted.
\qed
\end{defn}

Our main result of this section is as follows.
\begin{thm}[Superposition principle]
\label{Superposi}
Assume that Assumption \ref{assump2} is in force.
Assume further that the CFP equation \eqref{FP for nonlinear Zakai} has a unique regular solution for some $p>1$. Then there exists a filtered probability space $(\Omega,\cF,\hQ^0)$, supporting a $(d+1)$-dimensional Brownian motion $(B^1, Y)$ and
a continuous $\cF$-adapted $\sM_1(\hR^d)$ valued process $(\n_t)_{t\in[0,T]}$ satisfying SPDE \eqref{Zakai: state} with
  $\bP^Y_t=\hP^Y_t\circ (\n_t)^{-1}$, $\hQ^0$-a.s., for each $t\in[0,T]$.
\end{thm}

We shall follow the strategy that is similar to those in \cite{lacker2020, Qiao2022,  trevisan2014well}    to prove Theorem \ref{Superposi}. That is, we shall   transfer the CFP equation \eqref{FP for nonlinear Zakai} to a PDE on $\hR^\infty$ and then 
apply  the so-called {\it filtered martingale problem} proposed in \cite{kurtz1988} to obtain the weak solution to the nonlinear Zakai equation.

To begin with, let $\{\f_n: n\in\hN \}\subset \hC_c^{\infty}(\hR^d)$ be  
such that for any $\f\in \hC_c^{\infty}(\hR^d)$, one can find a subsequence $(n_k)_{k\in\hN}$ such that $(\f_{n_k},\nabla\f_{n_k}, \nabla^2\f_{n_k})_{k\in\hN}$ converges uniformly to $(\f,\nabla\f,\nabla^2\f )$. (see, e.g., \cite[Lemma 6.1]{lacker2020} for the existence of such countable set $\{\f_n\}$). Next, recall the spaces $(\sM(\hR^d), {\bf d})$ and $\hR^{\infty}_\sM$  defined in \S2, and  the identification mapping $\sT : \sM(\hR^d)\rightarrow \hR^{\infty}$ defined by (\ref{sT}). 
We note that, by using the mapping $\sT$ and its inverse $\sT^{-1}$  we can ``lift" the conditional Fokker-Planck equation \eqref{FP for nonlinear Zakai} into a conditional Fokker-Planck equation on $\hR^{\infty}$ as we briefly indicated in \S2. 
More precisely, for $t\in[0,T]$, let $\mathbf Q_t^Y$ be the
conditional probability law on $\sP(\hR^\infty_\sM)\subseteq \sP(\hR^\infty)$ induced by the mapping
$\sT$ under $\bP^Y_t$ as defined in \eqref{bPyt}. That is, $\bQ^Y_t=\bP^Y_t\circ \sT^{-1}=\hP^Y_t\circ (\n_t)^{-1}\circ\sT^{-1}\in\sP(\hR^\infty)$. Now for any $\bn\in\sM(\hR^d)$ and $ (\f_n)_{n\in\hN}\subset\hC^\infty_c(\hR^d)$, 
we denote $\bz=\sT(\bn)=\{\la \bn,\f_i\ra \}_{i=1}^\infty\in\hR^\infty$. Then, for any $k\in \hN$ and $f\in\hC(\hR^k) $, we derive the counterpart of conditional FPE on $\hR^\infty$ (\ref{R infinity FP}), which we rewrite here for ready reference: for $\hQ^0$-a.s.,
\bea
 \label{FRinfty}
\int_{\hR^{\infty}} f(\bz)(\mathbf Q^Y_t -\mathbf Q^Y_0)(d\bz)&=&\int_0^t\int_{\hR^{\infty }} 
\Big\{\sum_{i=1}^{k}\pa_if(\bz)\beta_i(s,Y_s,\bz) \Big\}\mathbf Q^Y_s(d\bz)ds\nonumber \\
&&+\int_0^t\int_{\hR^{\infty}}\frac{1}{2}\sum_{i,j=1}^k\pa_{ij}f(\bz)\alpha_{ij}(s,Y_s,\bz)\mathbf Q^Y_s(d\bz) ds\\
&&+\int_0^t\int_{\hR^{\infty }} 
\Big\{\sum_{i=1}^{k}\pa_if(\bz)\gamma_i(s,Y_s,\bz) \Big\}\mathbf Q^Y_s(d\bz)dY_s,  \q t\in[0,T], \nonumber\nonumber
\eea 
where, for $(s, Y_s, \bz)\in [0,T]\times \hR\times \hR^\infty$, 
{\small
\bea 
\label{alpha beta}
\left\{\ba{lll}
\alpha_{ij}(s, Y_s ,\bz) \neg &:=&\neg \big[\la \sT^{-1}(\bz),(\nabla\f_i(\cdot))^{\top}\rho(s, \cd,Y_s, \bm(\sT^{-1}(\bz)))+h(s, \cd, Y_s,\bm(\sT^{-1}(\bz)))\f_i(\cd)\ra\big]\\
 &&\neg\q\times \big[\la \sT^{-1}(\bz),(\nabla\f_j(\cdot))^{\top}\rho(s, \cd,Y_s, \bm(\sT^{-1}(\bz)))+h(s, \cd, Y_s, \bm(\sT^{-1}(\bz)))\f_j(\cd)\ra\big],\\ 
\beta_i(s,Y_s,\bz)  \neg&:=&\neg \la \sT^{-1}(\bz) ,L^Y_{s,\bz}\f_i (\cdot)\ra, \ss\\
 \gamma_i(s,Y_s,\bz)\neg  &:=&\neg  \la \sT^{-1}(\bz),(\nabla \f_i(\cd))^{\top}\rho(t, \cd,Y_{s}, \bm(\sT^{-1}(\bz)))+h(t, \cd, Y_{s}, \bm(\sT^{-1}(\bz)))\f_i(\cd)\ra,\\
L^Y_{s,\bz}\f_i\neg &:=&\neg \frac{1}{2}\si\si^{\top}(s, \cd , Y_s,\bm(\sT^{-1}(\bz) )):\nabla^2 \f_i(\cd)+(\nabla \f_i(\cd))^{\top}b(s, \cd , Y_s, \bm(\sT^{-1}(\bz) )).
\ea\right.  
\eea }

Clearly, for each $k\in\hN$, the coefficients $(\alpha_{ij})$, $(\beta_i)$, and $(\gamma_i)$ in \eqref{alpha beta} are  
Borel measurable functions defined on $[0,T]\times\hR\times \hR^\infty_{\sM}$. 
We now let $(\bP^Y_t)_{t\in[0,T]}$ be a given regular solution to the CFPE (\ref{FP for nonlinear Zakai}), and let $p>1$ be the corresponding index for the required integrability condition (\ref{integ}). 
We claim that, for the given $p>1$, the following integrability conditions hold $\hQ^0$-a.s.:   
\bea
\label{interg1}
 \int_0^T\int_{\hR^{\infty}}|\beta_i(t,Y_t,\bz)|^p\mathbf Q^Y_t(d\bz )dt&<\infty,\nonumber \\
 \int_0^T\int_{\hR^{\infty}}|\alpha_{ij}(t,Y_t,\bz) |^p \mathbf Q^Y_t(d\bz)dt&<\infty,\\
 \int_0^T\int_{\hR^{\infty}}|\gamma_i(s,Y_t,\bz)|^{2p} \mathbf Q^Y_t(d\bz)dt&<\infty. \nonumber
\eea
To prove (\ref{interg1}), first note that by definition (\ref{alpha beta}) and definition \eqref{bPyt} we have 
\beaa
&&\int_{\hR^{\infty}}|\beta_i(t,Y_t,\bz)|^p\mathbf Q^Y_t(d\bz )=\int_{\sM(\hR^d)}|\la \bn, L^Y_{t,\bz}\f_i\ra |^p\bP^Y_t(d\bn)\nonumber\\
&=&\int_{\sM(\hR^d) }\big|\la \bn, \frac{1}{2}\si\si^{\ts}(t, \cd , Y_t,\bm)):\nabla^2 \f_i(\cd)+(\nabla \f_i(\cd))^{\mathsf T}b(t, \cd , Y_t, \bm))\ra\big|^p\bP^Y_t(d\bn)  \\ 
&=&\int_{\sM(\hR^d)}\Big| \int_{\hR^d}\Big(\frac{1}{2}\si\si^{\ts}(t,x, Y_t,\bm)):\nabla^2 \f_i(x)+(\nabla \f_i(x))^{\mathsf T}b(t , x,Y_t, \bm)\Big)\bn(dx)\Big|^p\bP^Y_t(d\bn). \nonumber\\
&=&  \hE^{\hQ^0}\left[ \Big|\la \nu_t, \frac{1}{2}\si\si^{\ts}(t,x, Y_t,\m_t)):\nabla^2 \f_i(x)+(\nabla \f_i(x))^{\mathsf T}b(t , x,Y_t, \m_t)\ra \Big|^{p}\Big| \cF_t^Y  \right]\\
&\le&  \hE^{\hQ^0}\left[ \Big|\la |\nu_t|, \frac{1}{2}|\si\si^{\ts}(t,x, Y_t,\m_t)):\nabla^2 \f_i(x)|+|(\nabla \f_i(x))^{\mathsf T}b(t , x,Y_t, \m_t)|\ra \Big|^{p}\Big| \cF_t^Y  \right].
\eeaa

Now, by Assumption \ref{assump2} and the integrability condition   (\ref{integ}) we obtain
\bea
\label{beta}
&&\int_0^T\int_{\hR^{\infty}}|\beta_i(t,Y_t,\bz)|^p\mathbf Q^Y_t(d\bz )dt\nonumber\\
\neg&\neg\le\neg&\neg 2^{p-1}\int_0^T	 \int_{\sM(\hR^d)}\Big\{\big| \int_{\hR^d} \big[ \frac{1}{2}|\si\si^{\ts}(t,x, Y_t,\bm)):\nabla^2 \f_i(x)|\big] |\bn|(dx)\big|^p \\
&&\qq\qq+\big| \int_{\hR^d } |(\nabla \f_i(x))^{\top}b(t , x,Y_t,\bm)|\cdot|\bn|(dx) \big|^p\Big\}\bP^Y_t(d\bn)dt\nonumber\\
\neg&\neg\le\neg&\neg  2^{p-1}\| \f_i\|^p_{C_b^2 }\int_0^T\neg\neg\int_{\sM({\hR^d})}\neg\neg\Big[\|\frac{1}{2}\sigma\sigma^{\ts}(t,\cdot,Y_t,\bm)\|_{L^1(|\bn|)}^p+\|b(t , \cd,Y_t, \bm)\|^p_{L^1(|\bn|)}\Big]\bP^Y_t(d\bn)dt<\infty. \nonumber
\eea
Similarly, we have 
\beaa 
&& \int_{\hR^{\infty}}|\alpha_{ij}(t,Y_t,\bz) |^p \mathbf Q^Y_t(d\bz) =\int_{\sM(\hR^d)}\Big|\Big[\la \bn,(\nabla\psi_i(\cd))^{\mathsf T}\rho(s, \cd,Y_t, {\bm})+h(s, \cd, Y_t,{ \bm})\f_i(\cd)\ra\Big]\nonumber \nonumber\\
&&\qq\qq\qq\times \Big[\la \bn,(\nabla\f_j(\cd))^{\mathsf T}\rho(s, \cd,Y_t, {\bm})+h(s, \cd, Y_t,{ \bm})\f_j(\cd)\ra\Big]    \Big|^p \bp^Y_t(d\bn) \\
&=& \int_{\sM(\hR^d)}\Big|\int_{\hR^d} \Big[(\nabla\f_i(x))^{\mathsf T}\rho(t,x,Y_t, {\bm})+h(t, x,Y_t,{ \bm})\f_i(x)\Big] \bn(dx) \Big|^{p} \\
&& \times \Big|\int_{\hR^d} \Big[(\nabla\f_j(x))^{\mathsf T}\rho(t,x,Y_t, {\bm})+h(t, x,Y_t,{ \bm})\f_j(x)\Big] \bn(dx) \Big|^{p} \bp^Y_t(d\bn) \\
&\le & 2^{2(p-1)}\|\f_i\|^p_{C^1_b }\|\f_j\|^p_{C^1_b } \int_{\sM(\hR^d)} \Big| \|\rho(t,\cd,Y_t,\bm)\|^p_{L^1(|\bn|)} +\|h(t,\cd,Y_t,\bm)\|^p_{L^1(|\bn|)}\Big|^2\bP^Y_t(d\bn)\\
&\le& 2^{2p-1}\|\f_i\|^p_{C^1_b }\|\f_j\|^p_{C^1_b} \int_{\sM(\hR^d)}\Big[ \|\rho(t,\cd,Y_t,\bm)\|^{2p}_{L^1(|\bn|)} +\|h(t,\cd,Y_t,\bm)\|^{2p}_{L^1(|\bn|)}\Big]\bP^Y_t(d\bn).
\eeaa 
Now by Assumption \ref{assump2} and (\ref{integ}) again we have  
\bea
\label{alpha}
&&\int_0^T\int_{\hR^{\infty}}|\alpha_{ij}(t,Y_t,\bz) |^p \mathbf Q^Y_t(d\bz)dt\\
&\le& C_p(\f_i, \f_j)\int_0^T\int_{\sM(\hR^d)}\Big[ \|\rho(t,\cd,Y_t,\bm)\|^{2p}_{L^1(|\bn|)} +\|h(t,\cd,Y_t,\bm)\|^{2p}_{L^1(|\bn|)}\Big]\bP^Y_t(d\bn)dt<\infty,\nonumber
\eea
where $C_p(\f_i, \f_j):=2^{2p-1}\|\f_i\|^p_{C^1_b }\|\f_j\|^p_{C^1_b }$.
Similarly, by taking $i=j$ in $\alpha_{ij}$, we have,
\bea
 &&\int_{\hR^{\infty}}|\gamma_i(t,Y_y,\bz)|^{2p}\bQ_t^Y(d\bz) \nonumber\\
 &=&\int_{\sM(\hR^d)}|\la \bn,(\nabla \f_i(\cd))^{\mathsf T}\rho(t, \cd,Y_t, \bm)+h(t, \cd, Y_t, \bm)\f_i(\cd)\ra|^{2p}
\bP^Y_t(d\bn)\nonumber \\
&\le &C_p(\f_i)\int_{\sM(\hR^d)}\Big[ \|\rho(t,\cd,Y_t,\bm)\|^{2p}_{L^1(|\bn|)} +\|h(t,\cd,Y_t,\bm)\|^{2p}_{L^1(|\bn|)}\Big]\bP^Y_t(d\bn)<\infty \nonumber,
\eea
where $C_p(\f_i):=2^{2p-1}\|\f_i\|^{2p}_{C^1_b }.$
We  note that the integrability condition (\ref{interg1}) is a key element for us to carry out the  argument along the line of the so-called {\it filtered martingale problem} of  \cite{kurtz1988}. We are now ready to prove Theorem \ref{Superposi}. 

\ss
[{\it Proof of Theorem \ref{Superposi}.}] In view of our discussion on the connection between the CFP equation (\ref{FP for nonlinear Zakai}) and its $\hR^\infty$ counterpart (\ref{FRinfty}), we shall begin by assuming that there exists a $\sP(\hR^\infty)$-valued process $\bQ^Y_t$, $t\in[0,T]$, such that for any given $k\in\hN$, the SPDE (\ref{FRinfty})-(\ref{alpha beta}) hold for $\hQ^0$-a.s.. We are to show that the process 
$(\bQ^Y_t)_{t\in[0,T]}$ will lead us to a $\sP_1(\hC([0,T];\hR^\infty))$-valued random variable $\bQ^Y$ and a underlying $\{\cF^Y_t\}$-adapted process $V_t=(\bZ_t, Y_t)\in \hR^\infty\times \hR$, $t\in[0,T]$, such that $\n:=\sT^{-1}(\bZ)$ is a weak solution to the nonlinear Zakai equation
(\ref{Zakai: state}). That is, 
for any $\f\in \hC_c^{\infty}(\hR^d)$, the process $\lan \n_t, \f\ran$ satisfies (\ref{Zakai: state}), on a possibly extended probability space on which 
$Y$ is a Brownian motion.

To this end, let  $(\f_{n}) \subset \hC^\infty_c(\hR^d)$ be a sequence such that for any $\f\in \hC^\infty_c(\hR^d)$, there exists a subsequence $(\f_{n_\ell},\nabla\f_{n_\ell},\nabla^2\f_{n_\ell})$ converges uniformly to $(\f,\nabla\f,\nabla^2\f)$ as $\ell\to\infty$ (thanks to  \cite[Lemma 6.1]{lacker2020}). In what follows we shall only consider such a particular sequence as the basis for the $k$-cylindrical functions in (\ref{FP for nonlinear Zakai}). More precisely, we fix $k\in\hN$ and consider the $\hC^{k+1}_T$-valued process $(\bZ^{(k)},Y)$, in which $\bZ^{(k)}:=(\lan\sT^{-1}(\bZ), \f_1\ran, \cds, \lan\sT^{-1}(\bZ), \f_k\ran)$ is viewed as a $\hR^k$-valued {\it signal process} and $Y$ is the {\it observation process}.
For notational simplicity, in the rest of the argument we shall simply denote $\bZ=\bZ^{(k)}$ when the context is clear. We now briefly introduce the idea of {\it filtered martingale problem} for $\bZ$ given the observation  $\{\cF^Y_t\}_{t\in[0,T]}$, following the idea of \cite{kurtz1988}.

We begin by introducing the  generating operator $\mathcal A_0$ for the signal process $\bZ=\bZ^{(k)}$ 
associated to the equation (\ref{FRinfty}). For  $f\in\hC^{\infty}(\hR^{k+1})$, 
we define: for $\bz\in\hR^k$ and $\bar{y}\in\hR$, 
\beaa \label{FP generator}
\mathcal A_0[f](t,\bar{y}, \bz):=\frac{1}{2}\sum_{i,j=1}^k\pa_{ij}f(\bar{y}, \bz)\alpha_{ij}(t,\bar{y}, \bz)
+\sum_{i=1}^{k}\pa_if(\bar{y}, \bz)\beta_i(t,\bar{y}, \bz),
\eeaa 
where the derivatives $\partial_i$, $\partial_{ij}$ are taken with respect to the variable $\bz$, and
$\alpha_{ij} $ and $\beta$ are defined as in \eqref{alpha beta}.
 We then define the infinitesimal generator $\mathcal A$ for the filtered martingale problem associated with the coupled process $V:=(\bZ, Y)$ by
 \bea\label{operator A}
  \mathcal A[f](t,\bar{y}, \bz)= \mathcal A_0[f](t,\bar{y}, \bz)+ \frac{1}{2}\Delta_{y} [f](t,\bar{y}, \bz), \q (\bar{y}, \bz)\in\hR^{k+1}, 
   \eea
where $\Delta_{y}$ represents the Laplacian operator with respect to $\bar{y}$. Furthermore, we denote $\sD(\cA)$ (resp. $\sR(\cA)$) as the domain (resp. range) of the operator $\cA$. 
For notation convenience, let $V(t)=(\bZ_t,Y_t)\in\hR^k\times\hR$, for $k\in\hN$, and let $\mathbf Q_t$ denote the distribution of $V(t)$. A measurable process $V(t)$ is a solution of the martingale problem for $\mathcal A$ if there exits a filtration $\mathcal F_t$ such that 
\bea\label{martingale problem A}
f(V(t))-\int_0^t\mathcal A f(V(s))ds,
\eea
is an $\mathcal F_t$-martingale for each $f\in\cD(\cA).$ The martingale problem for the coupled process $V(t)=(\bZ_t,Y_t)$ with generator $\mathcal A$ is well-posed following from \cite[Theorem 6.4]{Qiao2022}. 
Then the martingale defined in \eqref{martingale problem A} implies that, for initial condition $p_0$,
\bea\label{martingale problem for A expectation}
\mathbf Q_t f=p_0 f +\int_0^t \mathbf Q_s \cA fds.
\eea 
Indeed, following \cite[Proposition 2.1]{kurtz1988}, one can show that given $p_0$, \eqref{martingale problem for A expectation} determines $\mathbf Q_t$ for $t\ge 0$. To show such uniqueness, we just need to prove that $\mathcal R(\lambda -\cA)$ is seperating for ecah $\lambda>0$, which is trivial here since $\cA$ is linear. Now, given that $V(t)=(\bZ_t,Y_t)$ is a solution of the martingale problem for $\mathcal A$, and $\mathbf Q^Y_t$ is the conditional distribution of $\bZ_t$ given $\mathcal F_t^Y$, we denote 
\beaa
\mathbf Q^Y_tf(\cdot, Y_t):=\int_{\hR^{k}} f(\bz,Y_t)\mathbf Q^Y_t (d\bz)=\hE[ f(\bZ_t,Y_t)|\mathcal F_t^Y],
\eeaa
for any $k\in\hN,$ and then 
\beaa
\int_{\hR^{k}} f(\bz,Y_t)\mathbf Q^Y_t (d\bz)-\int_0^t \int_{\hR^{k}}\mathcal A f(\bz,Y_{s}) \mathbf Q^Y_s(d\bz) ds,
\eeaa
is an $\{\mathcal F_t^Y\}$-martingale for each $f\in \mathcal D(\cA)$.  We then view $(\mathbf Q^Y, Y)$ as a process with sample paths in $\hC([0,T],\sP_1(\hR^k)\times \hR)$ for any $k\in \mathbb N$.  

Let us now recall the definition of a ``filtered martingale problem" for a $\hC^{k+1}_T$-valued pair $(\bZ, U)$ in the sense of  \cite{kurtz1988}. A $\hC([0,T],\sP_1(\hR^k)\times \hR)$-valued, continuous process $(\mathbf Q^U, U)$  is called a solution of the filtered martingale problem for $\mathcal A$ if $\mathbf Q ^U$ is $\hF^U$-adapted and
\bea 
\int_{\hR^{k}} f(\bz, U_t)\mathbf Q^U_t (d\bz)-\int_0^t \int_{\hR^{k}} \mathcal A f(\bz, U_s) \mathbf Q^U_s(d\bz) ds,\nonumber
\eea   
is an $\hF^U$-martingale for $f\in \sD(\cA)$.  The main difference is that our underlining SDE (i.e. the Zakai equation) are defined on $\hR^{\infty}\equiv \{\hR^k\}_{k\in\hN}$. Since $(\mathbf Q_t^U,U)$ is $\mathcal F_t^U$-adapted, for each $t$, there exists a Borel measurable function $H_t: \hC([0,T], \hR)\rightarrow \sP(\hR^k)$ such that 
\bea\label{H map}
\mathbf Q_t^U=H_t=H_t(U_{\cdot\wedge s})\quad a.s.,
\eea
for every $s>t$.  We then connect $(\mathbf Q^U,U)$ with the original problem $(\mathbf Q^Y,Y)$ by assuming that 
\bea\label{initial condition} 
\int_{\hR^{k}} f(\bz,U_0)\mathbf Q^U_0 (d\bz)=\int_{\hR^{k}} f(\bz,Y_0)\mathbf Q^Y_0 (d\bz),
\eea
for $f\in\mathcal B(\hR^k\times \hR)$, which means that $U_0$ and $Y_0$ have the same distribution. Next, we are to show the uniqueness of the filtered martingale problem. For process $U$ (and process $Y$ respectively), and $i=1,\cdots, m$,  we denote
\beaa\label{w function}
W^U_i(t)=\int_0^tg_i(U_s,s)ds=\int_0^tg_i(U_s, W_0(s))ds, \quad \text{for}\quad 0\le g_i\in \hC(\hR\times [0,\infty)) ,
\eeaa
where we make the convention $W_0(t)=t$.
Following  the same argument as in \cite{kurtz1988}[Lemma 3.1 \& Theorem 3.2], one can show that $(U, W_0^U,\cdots,W_m^U)$ and $(Y, W_0^Y,\cdots, W_m^Y)$ have the same one-dimensional distribution.  The idea is to construct a operator $\mathcal B$ acting on function class of the following type $f:\hR^k\times \hR\times \hR^{m+1}\rightarrow \hR
$,
\beaa
f(z,\tilde y,w_0,\cdots,w_m)=f(z,\tilde y)\Pi_{i=0}^mf_i(w_i),
\eeaa
where each $f_i$ denotes continuous functions vanishing at infinity. The operator $\mathcal B$ is defined below, 
\beaa
\mathcal B f(z,\tilde y, w_0,\cdots, w_m) =\mathcal A f(z,\tilde y)\Pi_{i=0}^mf_i(w_i)+f(z,y)\sum_{i=0}^mg_i(\tilde y,w_0)f_i'(w_i)\Pi_{i\neq j} f_j(w_j).
\eeaa
Following \cite{kurtz1988}, the existence and uniqueness of solution for the operator $\mathcal B$ follows from the existence and uniqueness of solution for the operator $\mathcal A$ \eqref{operator A}. We know that both $\alpha$ and $\beta$ are continuous Borel measurable maps and $\mathcal A$ is linear, for any $k\in\hN$, we know that $\mathcal R(\lambda-\mathcal A)$ is bounded-pointwise dense in $\mathcal B(\hR^k\times\hR)$ for each $\lambda >0$. Following \cite{kurtz1988}[Theorem 3.2], $\mathcal R(\lambda -\mathcal B)$ is separating for $\lambda>0$. Under such conditions, if $(\bZ_t,Y_t)$ is a solution of the $\hR^k\times\hR$ martingale problem for $\mathcal A$ and $(\mathbf Q_t^U,U)$ is a solution of the filtered margtingale problem for $\mathcal A$ with sample paths in $\hC([0,T],\sP_1(\hR^k)\times \hR)$ and satisfying \eqref{initial condition}, then there exists function $H_t$ satisfying \eqref{H map} and  $\mathbf Q_t^Y=H_t(Y)$  -$a.s.$ and $(\mathbf Q^Y, Y)$ has the same distribution as $(\mathbf Q^U,U)$. Since this is true for any $k\in\hN$,  we then conclude there exists a solution $\mathbf Q_t^Y\in\sP_1(C([0,T];\hR^{\infty}))$ of the filtered martingale problem associated to the system of SDEs,  
  \bea\label{SDE for Z}
d\mathbf Z_t^i&=&\la \sT^{-1}(\mathbf Z_t),L^{Y_{t}}_{t,\sT^{-1}(\mathbf Z_t)}\f_i(\cd)\ra dt+\la \sT^{-1}(\mathbf Z_t),(\nabla\f_i(\cd))^{\mathsf T}\rho(t,\cd,{Y_{t}}, {\bm(\sT^{-1}(\mathbf Z_t)) })\nonumber \\
&&+h(t,\cd, {Y_{t}}, \bm(\sT^{-1}(\mathbf Z_t)) )\f_i(\cd)\ra d Y_t,\quad \text{for}\quad i\in\hN.
\eea
Here we denote $\mathbf Z_t^i=\la \nu_t,\f_i\ra$ as the $i$-th projection of $\sT(\nu_t)$, and the corresponding marginal flow $\mathcal L(\mathbf Z_t|\mathcal F_t^Y)=\mathbf Q_t^Y$ denotes the solution of the conditional Fokker-Planck equation on $\hR^{\infty}$. 

\ms
Next we show that the solution $\mathbf Q_t^Y$ give rise to a weak solution of \eqref{SDE for Z}. For $\n_t=\sT^{-1}(Z_t)$ and $\m_t=\bm(\sT^{-1}(\mathbf Z_t))$, we have, for $i\in\hN$, 
\beaa
\la \n_t,\f_i\ra =\mathbf Z^i_t=\mathbf Z^i_0+\int_0^t \la \n_s,L^Y_{s,\m_s}\f_i\ra ds+\int_0^t\la \n_s,(\nabla\f_i)^{\mathsf T}\rho(s,\cd,Y_{s}, \m_s)+h(s,\cd, Y_{s}, \m_s)\f_i\ra d Y_s.
 \eeaa 
 Taking conditional expectation $\hE^{\hQ^0}[\cdot |\cF_t^Y ]$ on both sides above, noting that  all processes involved are $\hF^Y$-adapted, and applying Proposition \ref{prop: fokker plack}, we have 
 \bea \label{SDE for Z-II}
\mathbf Z^i_t&=&\mathbf Z^i_0+\int_0^t \hE^{\hQ^0}[ \la \n_s,L^Y_{s,\m_s}\f_i\ra |\cF_s^Y] ds\nonumber \\
&&+\int_0^t\hE^{\hQ^0}[ \la \n_s,(\nabla\f_i)^{\mathsf T}\rho(s,\cd,Y_{s}, \m_s)+h(s,\cd, Y_{s}, \m_s)\f_i\ra |\cF_s^Y] d Y_s \nonumber\\
&=& \mathbf Z^i_0+ \int_0^t \int_{\sM(\hR^d)} \la \bn,L^Y_{s,\bm}\f_i\ra \bP_s^Y(d\bn)  ds \nonumber\\
&&+\int_0^t  \int_{\sM(\hR^d)}\la \bn,(\nabla\f_i)^{\mathsf T}\rho(s,\cd,Y_{s}, \bm)+h(s,\cd, Y_{s}, \bm)\f_i\ra  \bP_s^Y(d\bn)d Y_s.
 \eea 
Now for any $\f\in \hC_c^{\infty}(\hR^d)$, by definition of the sequence $\{\f_n\}$, we can find a subsequence $(\f_{n_\ell})_{\ell\in\mathbb N}$ such that $(\f_{n_\ell},\nabla\f_{n_\ell},\nabla^2\f_{n_\ell})$ converges uniformly to $(\f,\nabla\f,\nabla^2\f)$ as $\ell\to\infty$. We claim that $\bZ_t:=\la\n_t,\f  \ra=\lim_{\ell\to\infty}  \la \n_t,\f_{n_\ell}\ra=\lim_{\ell\to\infty}\bZ^{n_\ell}_t$ a.s. for every $t\in[0,T]$, satisfies \eqref{Zakai: state}. Indeed, following the estimate in \eqref{beta}, we have
\bea
&& \Big| \int_0^t \int_{\sM(\hR^d)}\la \bn, L^Y_{s,\bm}\f_{n_\ell}\ra \bP_s^Y(d\bn) ds- \int_0^t\int_{\sM(\hR^d)}\la \bn, L^Y_{s,\bm} \f\ra\bP_s^Y(d\bn) ds  \Big| \nonumber \\
&\le & \| \f_{n_\ell}-\f\|_{C_b^2(\hR^d)}\int_0^t\neg\neg\int_{\sM({\hR^d})}\neg\neg\Big[\|\frac{1}{2}\sigma\sigma^{\ts}(s,\cdot,Y_s,\bm)\|_{L^1(|\bn|)}+\|b(s , \cd,Y_s, \bm)\|_{L^1(|\bn|)}\Big]\bP^Y_s(d\bn)ds.  \label{est L} \nonumber \\
\eea
The integral part of \eqref{est L}, for $t\in[0,T]$, is finite a.s. according to the integrability assumption for $\sigma,b,\rho,h$ in \eqref{integ}, by taking $\ell\rightarrow \infty$, we get \eqref{est L} $\rightarrow 0.$
Next, let $\Phi(t)$, $t\in [0,T]$, be defined as in \eqref{integ}. Define $\tau_N=\inf\{t>0, \Phi(t)>N \}$, then $\hE^{\hQ^0}[\Phi(T\wedge \tau_N)]\le N<\infty.$ For each $N$, applying It\^o isometry, we have 
\beaa
&&\hE^{\hQ^0}\Big[\Big( \int_0^{T\wedge \tau_N} \la \n_s,(\nabla\f_{n_\ell})^{\mathsf T}\rho(s,\cd,Y_{s}, \m_s)+h(s, \cd, Y_{s}, \m_s)\f_{n_\ell}\ra d Y_s\\
&&\qq\qq\qq- \int_0^{T\wedge \tau_N} \la \n_s,(\nabla\f)^{\mathsf T}\rho(s,\cd, Y_{s}, \m_s)+h(s,\cd, Y_{s}, \m_s)\f\ra d Y_s \Big)^2\Big]\\
&\le& \|\f_{n_\ell}-\f\|^2_{C_b^1(\hR^d)} \hE^{\hQ^0}\Big[ \int_0^{T\wedge\tau_N} \la|\n_s|, |\rho(s,\cd,Y_{s},\m_s)+|h(s,\cd,Y_{s},\m_s)|\ra^2 ds \Big] \\
&\le& \|\f_{n_\ell}-\f\|^2_{C_b^1(\hR^d)} \hE^{\hQ^0}\Big[ \int_0^{T} \hE^{\hQ^0}\Big[\mathbf{1}_{\{s\le \tau_N \}} \la|\n_s|, |\rho(s,\cd,Y_{s},\m_\ell)+|h(s,\cd,Y_{s},\m_s)|\ra^2 \Big| \cF_s^Y \Big]ds\Big] \\
&\le& \|\f_{n_\ell}-\f\|^2_{C_b^1(\hR^d)} \hE^{\hQ^0}\Big[ \int_0^{T\wedge\tau_N} \int_{\sM(\hR^d)}\Big(\int_{\hR^d}[|\rho(s,\cd,Y_{s},\bm)|+|h(s,\cd,Y_{s},\bm)| ]|\bn|(dx)\Big)^2 \bP_s^Y(\bn) ds \Big] \\
&\le&2 \|\f_{n_\ell}-\f\|^2_{C_b^1(\hR^d)} \hE^{\hQ^0}\Big[ \int_0^{T\wedge\tau_N} \int_{\sM(\hR^d)}
\Big(\|\rho(s,\cd,Y_{s},\bm)\|_{L^1(|\bn|)}^2+\|h(s,\cd,Y_{s},\bm)\|_{L^1(|\bn|)}^2\Big) \bP_s^Y(\bn) ds \Big]\\
&\le & 2 \|\f_{n_\ell}-\f\|^2_{C_b^1(\hR^d)} \hE^{\hQ^0} \Big[\Phi(T\wedge\tau_N)\Big] \le 2N \|\f_{n_\ell}-\f\|^2_{C_b^1(\hR^d)}  \rightarrow 0, \q \text{as} \q \ell\to \infty.
\eeaa
Therefore, let $\ell\to\infty$, we obtain,
\bea \label{SDE for Z-III}
\mathbf Z_t&=& \mathbf Z_0+ \int_0^t \int_{\sM(\hR^d)} \la \bn,L^Y_{s,\bm}\f\ra \bP_s^Y(d\bn)  ds \\
&&+\int_0^t  \int_{\sM(\hR^d)}\la \bn,(\nabla\f)^{\mathsf T}\rho(s,\cd,Y_{s}, \bm)+h(s,\cd, Y_{s}, \bm)\f\ra  \bP_s^Y(d\bn)d Y_s,\q t\in [0,T\wedge\tau_N].\nonumber
 \eea 
 Since $\tau_N$ is monotone increasing, by taking $N\rightarrow\infty$, we see that \eqref{SDE for Z-III} holds for all $t\in[0,T]$ a.s..
Finally, note that $\bZ_t=\la \nu_t,\f\ra$ and 
\beaa
&& \int_0^t\neg \int_{\sM(\hR^d)} \la \bn,L^Y_{s,\bm}\f\ra \bP_s^Y(d\bn) ds 
\neg +\neg\int_0^t\neg  \int_{\sM(\hR^d)}\la \bn,(\nabla\f)^{\mathsf T}\rho(s,\cd,Y_{s}, \bm)\neg+\neg h(s,\cd, Y_{s}, \bm)\f\ra  \bP_s^Y(d\bn)d Y_s\\
&&= \int_0^t \hE^{\hQ^0}[ \la \n_s,L^Y_{s,\m_s}\f\ra |\cF_s^Y] ds+ \int_0^t\hE^{\hQ^0}[ \la \n_s,(\nabla\f)^{\mathsf T}\rho(s,\cd,Y_{s}, \m_s)+h(s,\cd, Y_{s}, \m_s)\f\ra |\cF_s^Y] d Y_s\\
&&= \hE^{\hQ^0}\left[ \int_0^t  \la \n_s,L^Y_{s,\m_s}\f\ra ds+  \int_0^t\la \n_s,(\nabla\f)^{\mathsf T}\rho(s,\cd,Y_{s}, \m_s)+h(s,\cd, Y_{s}, \m_s)\f\ra  d Y_s\Big|\cF_t^Y\right]\\
&&= \int_0^t  \la \n_s,L^Y_{s,\m_s}\f\ra ds+  \int_0^t\la \n_s,(\nabla\f)^{\mathsf T}\rho(s,\cd,Y_{s}, \m_s)+h(s,\cd, Y_{s}, \m_s)\f\ra  d Y_s.
 \eeaa 
 Since $\f\in\hC^{\infty}_c(\hR^d)$ is arbitraty, we  proved \eqref{Zakai: state}, whence the theorem. 
\qed

\bibliographystyle{abbrv}
\bibliography{FM}

\begin{thebibliography}{10}

\bibitem{AGSbook}
L.~Ambrosio, N.~Gigli, and G.~Savar{\'e}.
\newblock {\em Gradient flows: in metric spaces and in the space of probability
  measures}.
\newblock Springer Science \& Business Media, 2005.

\bibitem{barbu2020nonlinear}
V.~Barbu and M.~R{\"o}ckner.
\newblock From nonlinear {F}okker--{P}lanck equations to solutions of
  distribution dependent {SDE}.
\newblock {\em The Annals of Probability}, 48(4):1902--1920, 2020.

\bibitem{bensoussan}
A.~Bensoussan.
\newblock {\em Stochastic control of partially observable systems}.
\newblock Cambridge University Press, 2004.

\bibitem{BLM}
R.~Buckdahn, J.~Li, and J.~Ma.
\newblock A mean-field stochastic control problem with partial observations.
\newblock {\em The Annals of Applied Probability}, 27(5):3201--3245, 2017.

\bibitem{buckdahn2021general}
R.~Buckdahn, J.~Li, and J.~Ma.
\newblock A general conditional {M}c{K}ean-{V}lasov stochastic differential
  equation.
\newblock {\em The Annals of Applied Probability.}, 33(3):2004--2023, 2023.

\bibitem{cardaliaguet2010notes}
P.~Cardaliaguet.
\newblock Notes on mean field games.
\newblock Technical report, 2010.

\bibitem{carmona2018probabilistic}
R.~Carmona, F.~Delarue, et~al.
\newblock {\em Probabilistic theory of mean field games with applications
  I-II}.
\newblock Springer, 2018.

\bibitem{coghi2019}
M.~Coghi and B.~Gess.
\newblock Stochastic nonlinear {F}okker--{P}lanck equations.
\newblock {\em Nonlinear Analysis}, 187:259--278, 2019.

\bibitem{crisan2014conditional}
D.~Crisan, T.~G. Kurtz, and Y.~Lee.
\newblock Conditional distributions, exchangeable particle systems, and
  stochastic partial differential equations.
\newblock In {\em Annales de l'IHP Probabilit{\'e}s et statistiques},
  volume~50, pages 3, 946--974, 2014.

\bibitem{hammersley2019weak}
W.~R. Hammersley, D.~Siska, and L.~Szpruch.
\newblock Weak existence and uniqueness for {M}c{K}ean-{V}lasov {SDE}s with
  common noise.
\newblock {\em The Annals of Probability}, 49:527--555, 2021.

\bibitem{kurtz1988}
T.~G. Kurtz and D.~L. Ocone.
\newblock Unique characterization of conditional distributions in nonlinear
  filtering.
\newblock {\em The Annals of Probability}, 16(1):80--107, 1988.

\bibitem{kurtzxiong}
T.~G. Kurtz and J.~Xiong.
\newblock Particle representations for a class of nonlinear {SPDE}s.
\newblock {\em Stochastic Processes and their Applications}, 83(1):103--126,
  1999.

\bibitem{lacker2018}
D.~Lacker.
\newblock Dense sets of joint distributions appearing in filtration
  enlargements, stochastic control, and causal optimal transport.
\newblock {\em arXiv preprint arXiv:1805.03185}, 2018.

\bibitem{lacker2020}
D.~Lacker, M.~Shkolnikov, and J.~Zhang.
\newblock Superposition and mimicking theorems for conditional
  {M}c{K}ean--{V}lasov equations.
\newblock {\em Journal of the European Mathematical Society}, 2022.

\bibitem{protter}
P.~E. Protter.
\newblock {\em Stochastic differential equations}.
\newblock Springer, 2005.

\bibitem{Qiao2022}
H.~Qiao.
\newblock Superposition principles for the {Z}akai equations and the
  {F}okker-{P}lanck equations on measure spaces.
\newblock {\em Bulletin des Sciences Mathématiques}, 174:103084, 2022.

\bibitem{ren2020space}
P.~Ren and F.-Y. Wang.
\newblock Space-distribution {PDE}s for path independent additive functionals
  of {M}c{K}ean--{V}lasov {SDE}s.
\newblock {\em Infinite Dimensional Analysis, Quantum Probability and Related
  Topics}, 23(03):2050018, 2020.

\bibitem{trevisan2014well}
D.~Trevisan.
\newblock {\em Well-posedness of diffusion processes in metric measure spaces}.
\newblock PhD thesis, PhD thesis, 2014.

\end{thebibliography}

\end{document}